\documentclass[10pt]{amsart}

\usepackage[hmargin=0.8in,twosideshift=0in,height=8.6in]{geometry}
\usepackage{amssymb,amsthm}
\usepackage{delarray,verbatim}
\usepackage{natbib}
\renewcommand{\cite}{\citet}

\usepackage{ifpdf}
\ifpdf 
\usepackage[pdftex]{graphicx}
\else
\usepackage[dvips]{graphicx}
\fi

\usepackage{ifpdf}
\usepackage{color}
\definecolor{webgreen}{rgb}{0,.5,0}
\definecolor{webbrown}{rgb}{.8,0,0}
\definecolor{emphcolor}{rgb}{0.95,0.95,0.95}

\usepackage{hyperref}
\hypersetup{%
          colorlinks=true,
          linkcolor=webbrown,
          filecolor=webbrown,
          citecolor=webgreen,
          breaklinks=true}
\ifpdf
\hypersetup{pdftex,
            pdfstartview=FitH, %%Fit, FitB, FitH
            bookmarksopen=true,
            bookmarksnumbered=true
}
\else
\hypersetup{dvips}
\fi

\usepackage{multicol}
\usepackage{multirow}

\renewcommand{\theequation}{\thesection.\arabic{equation}}
\numberwithin{equation}{section} 

\newtheorem{theorem}{Theorem}[section]
\newtheorem{lemma}{Lemma}[section]
\newtheorem{corollary}{Corollary}[section]
\newtheorem{remark}{Remark}[section]
\newtheorem{proposition}{Proposition}[section]
\newtheorem{definition}{Definition}[section]

\newcommand{\D}{\mathcal{D}}
\newcommand{\F}{\mathcal{F}}
\newcommand{\C}{\mathcal{C}}

\renewcommand{\P}{\mathbb{P}}
\newcommand{\E}{\mathbb{E}}

\newcommand{\eps}{\varepsilon}
\newcommand{\R}{\mathbb{R}}
\renewcommand{\S}{\mathcal{S}}

\renewcommand{\L}{\mathcal{L}}

\title[Analysis of the optimal exercise boundary]{Analysis of the optimal exercise boundary of \\ American options for jump diffusions}

\author{Erhan Bayraktar}
\address[Erhan Bayraktar]{Department of Mathematics, University of Michigan, Ann Arbor, MI 48109, USA; e-mail:erhan@umich.edu}
\thanks{This research is partially supported by the National Science Foundation.}
\thanks{We would like to thank Baojun Bian and Sijue
Wu for helpful discussions. We also would like to thank the Corresponding Editor Robert Pego, the anonymous Associate Editor and the two anonymous referees for their careful analysis of our paper. Their feedback helped us improve our paper.}

\author{Hao Xing }
\address[Hao Xing]{Department of Mathematics, University of Michigan, Ann Arbor, MI 48109, USA; e-mail:haoxing@umich.edu}

\begin{document}
\maketitle

\begin{abstract}
In this paper we show that the optimal exercise boundary / free
boundary of the American put option pricing problem for jump
diffusions is continuously differentiable (except at the
maturity). This differentiability result has been established by
Yang et al. (\textit{European Journal of Applied Mathematics}
17(1):95-127, 2006) in the case where the condition $r\geq q+
\lambda \int_{\R_+} \left(e^z-1\right) \nu(dz)$ is satisfied. We
extend the result to the case where the condition fails using a
unified approach that treats both cases simultaneously. We also
show that the boundary is infinitely differentiable under a
regularity assumption on the jump distribution.
\end{abstract}

\section{Introduction}

Let $(\Omega,\F,\P)$ be a complete probability space hosting a
Wiener process $W=\{W_t;t \geq 0\}$ and a Poisson random measure
$N$ on $\R_+ \times \R$ with the mean measure $\lambda dt\,
\nu(dz)$ (in which $\nu$ is a probability measure on $\R$)
independent of the Wiener process. Let $\mathbb{F}=\{\F_t\}_{t \in
[0,T]}$ be the (augmented) natural filtration of $W$ and $N$. We
will consider a Markov process $\mathbb{S}=\{S_t; t \geq 0\}$,
which follows the dynamics
\begin{equation}\label{eq:dyn}
dS_t=\mu S_{t-} dt+\widetilde{\sigma}(S_{t-}, t) S_{t-}
dW_t+S_{t-} \int_{\R} \left(e^z-1\right) N(dt,dz),
\end{equation}
as the stock price process. We will take $\mu \triangleq
r-q+\lambda-\lambda \xi$, in which
\begin{equation}\label{eq:xi}
\xi \triangleq \int_{\R}e^z \nu(dz) < \infty,
\end{equation}
as a standing assumption. We impose this condition on $\xi$ so
that the discounted stock prices are martingales. The constant
$r\geq 0$ is the interest rate, $q\geq 0$ is the dividend. The
volatility $\widetilde{\sigma}(S, t)$ is assumed to be
continuously differentiable in both $S$ and $t$. Moreover, there
are positive constants $\delta$ and $\Delta$ such that
\begin{equation}\label{eq:ass-tildesigma}
 0<\delta \leq \widetilde{\sigma}(S,t) \leq \Delta, \quad \text{ for all
 } S, t \geq 0.
\end{equation}
We should note that at the time of a jump the stock price moves
from $S_{t-}$ to $S_{t-} e^Z$ in which $Z$ is a random variable
whose distribution is given by $\nu$. When $Z<0$ the stock price
jumps down, when $Z>0$ the stock price jumps up. In the classical
Merton jump diffusion model, $Z$ is a Gaussian random variable.

In this framework, we will study the American put option pricing
problem. The value function of the American put option is defined
by
\begin{equation}\label{eq:opt-stop}
V(S,t) \triangleq \sup_{\tau \in \S_{0,T-t}}\E\{e^{-r \tau}(K-S_{\tau}
)^+ | \,S_0= S\},
\end{equation}
in which $\S_{0,T-t}$ is the set of stopping times (with respect
to the filtration $\mathbb{F}$) taking values in $[0,T-t]$. The
value function $V$ is the classical solution of a free boundary
problem (see Proposition~\ref{prop:classical-soln}). The main goal
of this paper is to analyze the regularity of the free boundary.
We will show that the free boundary is $C^1$ except at the
maturity $T$, and $C^{\infty}$ with an appropriate regularity
assumption on the jump distribution $\nu$. For notational
simplicity we will first change variables and transform the value
function $V$ into $u$ and its free boundary $s$ into $b$ (see
(\ref{eq:change-var})) and state our results in terms of $u$ and
$b$.

While the continuity of the free boundary of the American put option
in jump models has been studied extensively, for example, by
\cite{pham}, \cite{yang} and \cite{lamberton}, the 
differentiability of the free boundary was left as an open problem.
Even when the geometric
Brownian motion is the underlying process the differentiability is difficult to
establish (see the discussion on page 172 of \cite{peskir}) and
has only recently been fully analyzed by \cite{chen-chadam}. In
the jump diffusion case, \cite{yang} proved that the free boundary
is continuously differentiable before the maturity when the
parameters satisfy
\begin{equation} \label{eq:yang-cond}
r\geq q+ \lambda \int_{\R_+} \left(e^z-1\right) \nu(dz).
\end{equation}
When the condition \eqref{eq:yang-cond} is violated, the free
boundary of the American option for jump diffusions exhibits a
discontinuity at the maturity (see Theorem 5.3 in \cite{yang} and
equation (\ref{eq:b(0)}) in this paper). This behavior of the free
boundary was also observed by \cite{leven} and \cite{lamberton} in
the exponential L\'{e}vy models. The purpose of our paper is to
extend the regularity results of the free boundary to the case
where (\ref{eq:yang-cond}) is not satisfied. We will see that the
boundary is differentiable even when
\eqref{eq:yang-cond} is violated.

There are two critical points in showing the differentiability
properties without the condition \eqref{eq:yang-cond}: 1) to show
the H\"{o}lder continuity of the free boundary, 2) to show that
$\partial^2_S V(S,t)$ is strictly larger than 0 when the point
$(S,t)$ is close to the free boundary in the continuation region.
We achieve these two results in Theorem~\ref{theorem:holder-cont}
and Corollary~\ref{cor:u_xx_pos} respectively. Combining these two
properties and a generalization of the result in \cite{Cannon}
(see Lemma~\ref{lemma:u_xt cont}), we upgrade the regularity of
the free boundary from H\"{o}lder continuity to continuous
differentiability in Theorem~\ref{theorem:cont-diff}. Then we
analyze the higher order regularity of the free boundary making
use of a technique \cite{schaeffer} used for the free boundary of
a one dimensional Stefan problem on a bounded domain.

In order to show that the free boundary is continuously differentiable, it is essential that the value function $V(S,t)$
is the unique classical solution of the free boundary problem and has a continuous second derivative (see \eqref{eq:b'(t)}). 
In the
jump diffusion models, this has been shown by \cite{pham} under
 condition \eqref{eq:yang-cond}.
This condition was removed in \cite{yang} and also in
\cite{bayraktar-finite-horizon}. Moreover, continuous differentiability of the free boundary requires the continuity of the cross derivatives of the value function.
In the L\'{e}vy models with
infinite activity jumps, the value function is not expected to be
a classical solution in general. Yet in the literature different
notions of generalized solutions were explored. For example,
\cite{pham-vis} showed that the value function is a viscosity solution, \cite{achdou} showed
that the value function is the solution in the Sobolev sense and
\cite{lamberton} proved that the value function is the solution in
the distribution sense. Moreover, the smooth-fit property (see
\eqref{eq:V_3}) is also necessary in our analysis (see
Theorem~\ref{theorem:cont-diff} and equation \eqref{eq:u_3_t}).
While this property may not hold for general pay-off functions
(see \cite{peskir-angle}), it has been shown to hold for the put
option pay-off in \cite{zhang}, \cite{pham} and
\cite{bayraktar-finite-horizon} in the jump diffusion models. The
analysis in this paper also applies to the pay-off functions which
are continuously differentiable, bounded, convex on $[0,+\infty)$
and equal to zero in $[K,+\infty)$. In fact, the
singularity at the strike of the put option pay-off is the source of the
technical difficulties. Therefore, we will focus on the put option
pay-off in this paper and leave the investigation of the boundary behavior for general
pay-off functions to future work.

The rest of the paper is organized as follows:
 In Section 2, after changing variables we will collect several useful properties  of
the function $u$, which will be crucial in establishing our main
results in the next three sections about the regularity of its
free boundary. In Section 3, we will introduce an auxiliary
function and use it to show that the the free boundary is
H\"{o}lder continuous. In Section 4, we will prove the continuous
differentiability of the free boundary. In Section 5, we will
upgrade the regularity of the boundary curve and show that it is
infinitely differentiable under an appropriate regularity
assumption on the jump distribution. Finally, in Section 6, we
will show that the approximation free boundaries, constructed in
\cite{bayraktar-finite-horizon}, have the similar regular
properties with the original free boundary. Proofs of some
auxiliary results are presented in the Appendix.

\begin{figure}[t]\label{figure:theorem-diag}
\begin{minipage}{\textwidth}
\begin{center}
\caption{Our results and the relationships among them. } $A
\rightarrow B$ means that statement A is used in the proof of
statement B.
\includegraphics[width=6in, height=3in]{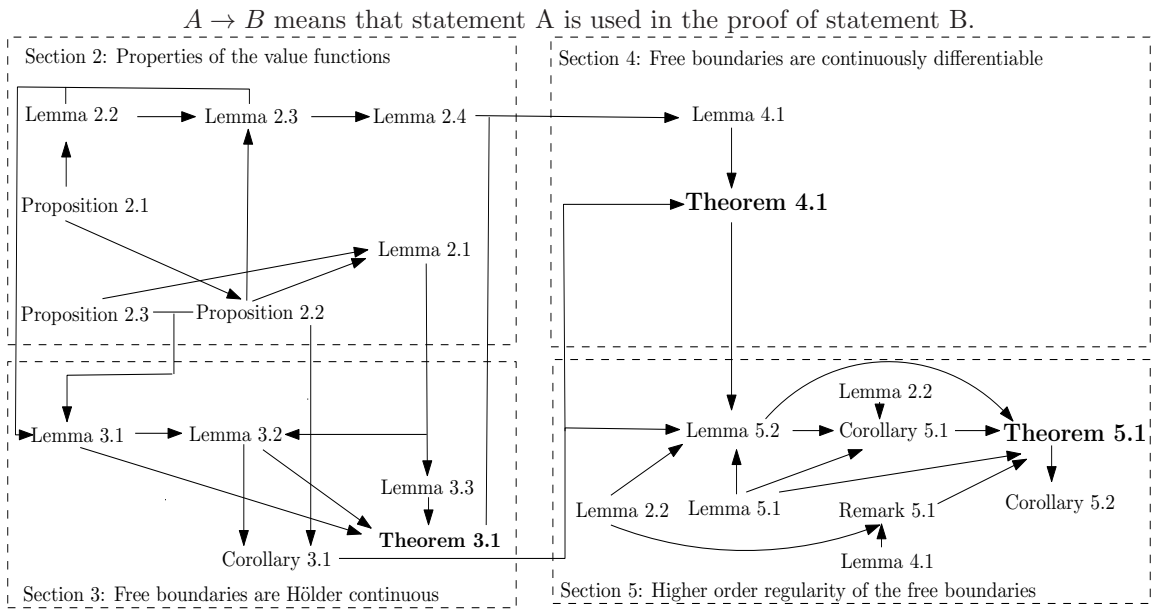}
\end{center}
\end{minipage}
\end{figure}

Our main results are Theorems~\ref{theorem:holder-cont},
\ref{theorem:cont-diff} and \ref{theorem:higher-reg-b}. In Figure
\ref{figure:theorem-diag} we show the logical flow of the paper,
i.e. we show how several results proved in the paper are related
to each other.

\section{Properties of the value function}

The value function $V(S,t)$ of the American put option for jump
diffusions solves a free boundary problem with the free boundary
$s(t)$. In particular, Theorem 4.2 of \cite{yang} and Theorem 3.1
of \cite{bayraktar-finite-horizon} state the following:

\begin{proposition}\label{prop:classical-soln}
$V(S,t)$ is the unique classical solution of the following boundary value problem:
\begin{eqnarray}
 &&\frac{\partial V}{\partial t} + \frac12 \widetilde{\sigma}(S,t)^2 S^2 \frac{\partial ^2
 V}{\partial S^2} + \mu S\frac{\partial V}{\partial S} - (r+\lambda)V +
 \lambda\int_{\R}V(S e^z,t)\nu(dz) =0, \, S>s(t), \label{eq:V_1}\\
 &&V(s(t),t)=K-s(t), \quad t\in [0,T), \label{eq:V_2}\\
 && V(S,T) = (K-S)^+, \quad S\geq s(T). \label{eq:V_4}
\end{eqnarray}
Moreover, the smooth fit property is satisfied, i.e.
\begin{equation}
 \frac{\partial}{\partial S}V(s(t), t) =-1, \quad t\in[0,T). \label{eq:V_3}
\end{equation}
In the region $\{(S,t) : S<s(t), t\in[0,T)\}$, $V(S,t)$ also
satisfies the following inequality:
\begin{equation}\label{eq:V_5}
 \frac{\partial V}{\partial t} + \frac12 \widetilde{\sigma}(S,t)^2 S^2 \frac{\partial ^2
 V}{\partial S^2} + \mu S\frac{\partial V}{\partial S} - (r+\lambda)V +
 \lambda\int_{\R}V(S e^z,t)\nu(dz) \leq0.
\end{equation}
\end{proposition}

In the following, let us first change the variables to state
(\ref{eq:V_1})-(\ref{eq:V_5}) in a more convenient form:
\begin{equation} \label{eq:change-var}
 x=\log(S), \quad u(x,t)=V\left(S,T-t\right), \quad b(t)=\log\left(s\left(T-t\right)\right)
\quad \text{  and  } \quad \sigma(x,t)= \widetilde{\sigma}(S,t).
\end{equation}
It is clear from the assumptions of $\widetilde{\sigma}(S,t)$ that
\begin{equation}\label{eq:ass-sigma}
\begin{split}
 & \sigma \text{ is continuously differentiable in both variables and} \\
 & \text{there are positive constants } \delta \text{ and } \Delta \text{ such that } 0< \delta <
 \sigma(x,t) <\Delta \text{ for all } (x,t)\in \R\times[0,T].
\end{split}
\end{equation}
While the first part of \eqref{eq:ass-sigma} will be used in 
\eqref{eq:def-h} and Lemma~\ref{lemma:u_xt cont}, the second part, which makes sure that the differential operators involved are uniformly parabolic,
will be necessary for Lemma~\ref{lemma:ut-bounded},
Corollary~\ref{cor:ut-holder} and
Theorem~\ref{theorem:higher-reg-b}. For the simplicity of the
notation, we will omit the variables of $\sigma$ in the sequel.
In terms of the new variables introduced in \eqref{eq:change-var}, (\ref{eq:V_1}) - (\ref{eq:V_5})
reduce to the uniformly parabolic boundary value problem
\begin{eqnarray}
 &&\L u \triangleq \frac{\partial u}{\partial t} -
 \frac12 \sigma^2 \frac{\partial^2 u}{\partial x^2} -
 \left(\mu - \frac12 \sigma^2\right)\frac{\partial u}{\partial x}
 + (r+\lambda)u - \lambda \int_{\R}u(x+z,t) \nu(dz)=0,  \quad x>b(t), \label{eq:u_1}\\
 && u(b(t),t)= K-e^{b(t)} ,\quad t\in (0,T],
 \label{eq:u_2}\\
 && u(x,0)=(K-e^x)^+, \quad x\geq b(0), \label{eq:u_4}\\
 && \frac{\partial}{\partial x} u(b(t),t) = -e^{b(t)},
 \label{eq:u_3} \\
 && \L u(x,t) \geq 0, \quad x<b(t), \, t\in(0,T]. \label{eq:u_5}
\end{eqnarray}
Let us define the continuation region $\mathcal{C}$ and the
stopping region $\mathcal{D}$ as follows
\begin{eqnarray*}
\mathcal{C} \triangleq \{(x,t)\,| \,b(t)<x<+\infty, 0<t\leq T\},
\quad  \mathcal{D} \triangleq \{(x,t)\,|\, -\infty<x\leq b(t),
0<t\leq T\}.
\end{eqnarray*}
From Proposition \ref{prop:classical-soln}, it is clear that the
boundary value problem (\ref{eq:u_1}) - (\ref{eq:u_4}) has a
unique classical solution $u(x,t)$ in $\mathcal{C}$.

\begin{remark}\label{remark:para-eq-u}
The integral term in (\ref{eq:u_1}) can also be considered as a
driving term, then the integro-differential equation
(\ref{eq:u_1}) can be viewed as the following parabolic
differential equation with a driving term $f(x,t)=\lambda
\int_{\R} u(x+z,t)\nu(dz)$:
\begin{equation} \label{eq:u_1_para}
\frac{\partial u}{\partial t} -
 \frac12 \sigma^2 \frac{\partial^2 u}{\partial x^2} -
 \left(\mu-\frac12 \sigma^2\right)\frac{\partial u}{\partial x}
 + (r+\lambda)u = f(x,t).
\end{equation}
 This point of view will be useful in the proof
of some results in later sections.
\end{remark}

In this section, we will study the properties of $u$ in both the
continuation and the stopping regions. Let us start from the
following proposition from \cite{yang}. It shows that the time
derivative of $u$ is continuously differentiable across the free
boundary.

\begin{proposition}\label{prop:u_t}
$\partial_t u(x,t)$ is a continuous function in $\R\times (0,T]$.
In particular, for any $t\in(0,T]$,
\begin{eqnarray}\label{eq:cont-at-b-us}
 \lim_{x\downarrow b(t)} \frac{\partial}{\partial t} u(x,t) &=&0
 \label{eq:u_t}.
 \end{eqnarray}
\end{proposition}
\begin{proof}
The proof is given in Theorem 5.1 in \cite{yang}, which summarized
Lemmas 2.8 and 2.11 in the same paper and used a special case of
Lemma 4.1 in page 239 of \cite{friedman-2}.
\end{proof}

Moreover, we will show in the following that $t\rightarrow u(x,t)$
is strictly increasing function in the continuation region.
\begin{proposition}\label{prop:u_t positive}
\begin{eqnarray}
 &&\frac{\partial u}{\partial t}(x,t) >0, \quad (x,t)\in
 \mathcal{C}.
 \label{eq:u_t_positive}
\end{eqnarray}
\end{proposition}

\begin{proof}
 The inequality (\ref{eq:u_t_positive}) is proved in Proposition
 4.1 in \cite{yang} using the Maximum Principle for the
 integro-differential equations, which can be found in Theorem 2.7
 in Chapter 2 of \cite{garr-mena-1}. However, it can be proved using the ordinary Maximum Principle for
 parabolic differential equations (see Corollary 7.4 in \cite{lieberman}). We know that $w=\frac{\partial u}{\partial
 t}$ satisfies the following equation in $\mathcal{C}$,
 \begin{eqnarray}
  \L_{\D}w &=& \lambda\int_{\R} w(x+z,t) \nu(dz), \label{eq:L-Dw}\\
 \end{eqnarray}
 Since $w=\partial_t u\geq 0$ in $\R\times(0,T)$, (\ref{eq:L-Dw}) implies that $\L_{\D} w \geq 0$. If there
 is a point $(x_0,t_0) \in \mathcal{C}$ such that $w(x_0,t_0)=0$ (i.e. $w$ achieves its non-positive minimum at $(x_0,t_0)$ ), it follows from
 the strong Maximum Principle that $w(x,t)=0$ in $\mathcal{C} \cap \R\times\{0<t\leq t_0\}
 $. Together with the fact that $w(x,t)=0$ in $\mathcal{D}$, we have that $w(x,t)=0$ in $\R\times \{0<t\leq t_0\}$.
 As a result, from
 \[
 u(x_0,t_0)-u(x_0,0) = \int_0^{t_0} w(x_0,s) ds =0,
 \]
 we obtain $u(x_0,t_0)= (K-e^{x_0})^+$. This contradicts with the definition of the free boundary
 $b(t)$, because $b(t_0)= \max\{x \in \mathbb{R}:u(x,t_0)=(K-e^{x})^+\}$ and
 $x_0>b(t_0)$.
\end{proof}

Combining Propositions~{\ref{prop:u_t}} and {\ref{prop:u_t
positive}} with the Hopf's Lemma for parabolic
integro-differential equations (see Theorem 2.8 in page 78 of
\cite{garr-mena-1}), we obtain that the free boundary is strictly
decreasing.

\begin{lemma}\label{lemma:b_strictly}
The function $t \to b(t)$ is strictly decreasing for $t\in(0,T]$.
\end{lemma}

\begin{proof}The proof is given in Theorem 5.4 in \cite{yang}.
\end{proof}

In order to investigate the regularity of the free boundary in the
later sections, we need more properties of $u$, which we will
develop in the following three lemmas. Since the results of these
lemmas are intuitive but proofs are technical, we will list the
proofs of these lemmas in the Appendix \ref{app:A.1}.

It is well
known that $S\rightarrow V(S, t)$ is uniformly Lipschitz in $\R_+$
and $t\rightarrow V(S, t)$ is uniformly semi-H\"{o}lder continuous
in $[0,T]$ (see \cite{pham}). The following
lemma shows the same properties also holds for $u(x,t)$, the
function that we obtained after the change of variables in
(\ref{eq:change-var}). (The globally Lipschitz continuity with
respect to $x$ is not a priori clear and one needs to check
whether $\partial_x u(x,t)$ is bounded.)

\begin{lemma}\label{lemma:lip-holder}
Let $u(x,t)$ be the solution of equation (\ref{eq:u_1}) -
(\ref{eq:u_4}), then we have
\begin{eqnarray}
&&|u(x,t)-u(y,t)|\leq C|x-y|, \quad x, y \in \R, t\in
[0,T],\label{eq:lip-u}\\
&&|u(x,t)- u(x,s)| \leq D|t-s|^{\frac12}, \quad x\in \R, 0\leq t,
s\leq T, \label{eq:holder-u}
\end{eqnarray}
where $C$ and $D$ are positive constants independent of $x$ and
$t$.
\end{lemma}
\begin{proof}
See Appendix \ref{app:A.1}.
\end{proof}

In the rest of this section, we will investigate the boundness of
$\partial_t u(x,t)$ and its behavior when $x\rightarrow +\infty$.
These two properties will be useful to show several results in
Sections~\ref{sec:cont-diff} and \ref{sec-inf-dif}  (see e.g.
(\ref{eq:w_1}), proof of Lemma~\ref{lemma:u_xt cont} and
Remark~\ref{remark:uniqueness}). Let us first recall the
definition of the H\"{o}lder spaces on page 7 of \cite{lad}.
\begin{definition}\label{def:holder space}
Let $\Omega$ be a domain in $\R$, $Q_T=\Omega\times (0,T)$. We
denote $\overline{Q_T}$ the closure of $Q_T$. For any positive
nonintegral real number $\alpha$,
$H^{\alpha,\alpha/2}\left(\overline{Q_T}\right)$ is the Banach
space of functions $v(x,t)$ that are continuous in
$\overline{Q_T}$, together with continuous derivatives of the form
$\partial_t^r\partial_x^s v$ for $2r+s<\alpha$, and have a finite
norm
\[
 ||v||^{(\alpha)} = |v|_x^{(\alpha)} + |v|_t^{(\alpha/2)} + \sum_{2r+s\leq [\alpha]} ||\partial_t^r\partial_x^s v||^{(0)},
\]
in which
\begin{eqnarray*}
||v||^{(0)} &=& max_{Q_T} |v|, \\
|v|_x^{(\alpha)} &=& \sum_{2r+s=[\alpha]}
<\partial_t^r\partial_x^s v>_x^{(\alpha-[\alpha])}, \quad
|v|_t^{(\alpha/2)} = \sum_{\alpha-2<2r+s<\alpha}
<\partial_t^r\partial_x^s v>_t^{(\frac{\alpha-2r-s}{2})};\\
<v>_x^{(\beta)} &=& \sup_{\begin{array}{cc}(x,t), (x',t)\in
\overline{Q_T}\\ |x-x'|\leq \rho_0\end{array}}
\frac{|v(x,t)-v(x',t)|}{|x-x'|^{\beta}}, \quad 0<\beta<1,\\
<v>_t^{(\beta)} &=& \sup_{\begin{array}{cc}(x,t), (x,t')\in
\overline{Q_T}\\ |t-t'|\leq \rho_0\end{array}}
\frac{|v(x,t)-v(x,t')|}{|t-t'|^{\beta}}, \quad 0<\beta<1,\\
\end{eqnarray*}
where $\rho_0$ is a positive constant.

On the other hand, $H^{\alpha}\left(\overline{\Omega}\right)$ is
the Banach space whose elements are continuous functions $f(y)$ on
$\overline{\Omega}$ that have continuous derivatives up to order
$[\alpha]$ and the following norm finite
\[
 || f ||^{(\alpha)} =  \sum_{j\leq [\alpha]} \left\| d_y^j f \right\|^{(0)}  + \left | d_y^{[\alpha]}f\right|^{(\alpha-[\alpha])} ,
\]
in which
\[
 |f|^{(\beta)}=  \sup_{y,y'\in \overline{\Omega}, |y-y'|\leq \rho_0} \frac{|f(y)-f(y')|}{|y-y'|^{\beta}}.
\]
Here $d_y^j f$ is the $j$th derivative of $f$. These H\"{o}lder norms depend on $\rho_0$, but for different $\rho_0>0$, the corresponding H\"{o}lder norms are equivalent hence their dependence on $\rho_0$ will not be noted in the sequel.
\end{definition}

Using the H\"{o}lder spaces and regularity results for parabolic
equations, we have the following result.
\begin{lemma}\label{lemma:ut-bounded}
For any $\epsilon>0$, $\partial_t u(x,t)$ is uniformly bounded on
$\R\times[\epsilon, T]$.
\end{lemma}
\begin{proof}
See Appendix \ref{app:A.1}.
\end{proof}

\begin{remark}\label{remark:condition_c1}
{\bf{(i)}} In the statement of Lemma \ref{lemma:ut-bounded}, $t=0$
cannot be included, i.e., $\lim_{t\rightarrow 0}\partial_t u(x,t)$
is not uniformly bounded in $x\in \R$, because  $\partial_t
u=\frac12 \sigma^2\partial^2_{x} u + \left(\mu-\frac12
\sigma^2\right)\partial_x u - (r+\lambda)u + \lambda\int_{\R}
u(x+z,t)\nu(dz)$ and $\lim_{t \rightarrow 0}
\partial^2_{x}u(x,t)$ is not bounded as a result of non-smoothness
of the initial value at $x=\log{K}$.
\end{remark}

In the following, we will use the previous lemma to  analyze the
behavior of $\partial_t u(x,t)$ as $x\rightarrow +\infty$.
\begin{lemma}\label{lemma:u_t_inf}
\[
 \lim_{x\rightarrow +\infty} \partial_t u(x,t) =
 0, \quad t\in (0,T].
\]
\end{lemma}
\begin{proof}
See Appendix \ref{app:A.1}.
\end{proof}

\begin{remark}\label{remark:u_c1}
Given the result in Lemma~\ref{lemma:ut-bounded}, it is clear from
the differential equation \eqref{eq:u_1_para} that $\partial^2_x
u$ is uniformly bounded in $\R\times[\epsilon, T]$, since
$\partial_x u$ is uniformly bounded (see
Lemma~\ref{lemma:lip-holder}). Combining with semi-H\"{o}lder
continuity of $u(x,\cdot)$ in Lemma~\ref{lemma:lip-holder}, Lemma
3.1 in page 78 of \cite{lad} now tells us that $\partial_x u(x,
\cdot)\in H^{1/2}([\epsilon, T])$. Therefore, combining with the
smooth fit property and Proposition~\ref{prop:u_t}, we have
\[
 u\in C^1 \left(\R\times (0, T]\right).
\]
\end{remark}

In the following three sections we will use the properties of the
value function we have shown in this section to investigate the
regularity of the free boundary $b(t)$.

\section{The free boundary is H\"{o}lder continuous}
\subsection{An auxiliary function}
Before we begin to analyze the regularity of the free boundary,
let us introduce the following important auxiliary function, which
was also used in \cite{lamberton} to prove the continuity of the
free boundary in an exponential L\'{e}vy model:
\begin{equation}\label{eq:def_J}
 J(x,t) \triangleq qe^x - rK + \lambda \int_{\R} \left[u(x+z, t) + e^{x+z}
 -K\right]\nu(dz), \quad x\in\R, \, t\in [0,T].
\end{equation}
As a result of the assumption \eqref{eq:xi}, $J<\infty$. Moreover, $J$ is closely related to the behavior of the value
function $u$ in the stopping region, since one can check that
\begin{eqnarray}
 \L u(x,t) &=& - J(x,t), \quad \text{ for } x<b(t),
 \, t\in(0,T], \label{eq:lu=j}\\
\L g(x) &=& \L u(x,0) = -\left[qe^x - rK +\lambda \int_{\R}
\left(e^{x+z} - K\right)^+ \nu(dz)\right] = - J(x,0) , \quad
\text{for } x < \log{K}, \label{eq:lg}
\end{eqnarray}
in which $g(x)\triangleq \left(K-e^x\right)^+$. As we shall see in
the rest of this section, the function $J(x,0)$ is of special
importance. We rename it as $J_0(x)$, i.e.,
\begin{equation}\label{eq:def-J0}
 J_0(x) \triangleq qe^x - rK +\lambda \int_{\R} \left(e^{x+z} - K\right)^+
 \nu(dz).
\end{equation}

Let us analyze the properties of $J$.
\begin{lemma}\label{lemma:j}
\begin{enumerate}
 \item $J(x,t) \geq -r K$, $\lim_{x\downarrow -\infty} J(x,t) =
 -rK$  and $\lim_{x\uparrow +\infty} J(x,t) = +\infty$,\\
 \item $J(x,t)\in C^1 \left(\R\times (0,T]\right)\cap C\left(\R\times[0,T]\right)$ ,\\
 \item
 The functions $x \rightarrow J(x,t)$ and $t\rightarrow J(x,t)$ are non-decreasing. If either either $q>0$ or
 \begin{eqnarray}
 \nu\left((M,+\infty)\right) >0, \quad \text{ for any} \quad M>0;
  \label{eq:condition-nu}
 \end{eqnarray}
 then $x \rightarrow J(x,t)$ is a strictly increasing function.
 On the other hand, if
 \begin{equation}\label{eq:condition-q}
 v((0,\infty))>0
 \end{equation}
  \begin{equation}\label{eq:jt}
    \partial_t J(x,t)  >0, \quad x\geq b(t), t\in(0,T].
  \end{equation}
\end{enumerate}
\end{lemma}
\begin{proof}
{\bf{(i)}} The first statement follows from $u(x+z,t) \geq
(K-e^{x+z})^+ \geq K-e^{x+z}$. The two limit statements follow
from the Bounded Convergence Theorem. \\{\bf{(ii)}}  The
continuity of $u(x,t)$ on $\R\times[0,T]$ implies that $J$ is
continuous on the same region. For the differentiability, since
$\partial_x u$ and $\partial_t u$ are uniformly bounded in
$\R\times[\epsilon, T]$ for any $\epsilon
>0$ (see Lemmas \ref{lemma:lip-holder} and
\ref{lemma:ut-bounded}), the Bounded Convergence Theorem gives us
\begin{equation}\label{eq:addt1}
\begin{split}
 &\frac{\partial}{\partial x}J(x,t) = q e^x + \lambda \int_{\R}
 \left[\frac{\partial}{\partial x} u(x+z, t) + e^{x+z}\right]
 \nu(dz) < +\infty,\\
 & \frac{\partial}{\partial t}J(x,t) = \lambda \int_{\R}
 \frac{\partial}{\partial t}u(x+z,t)\nu(dz) < +\infty.
\end{split}
\end{equation}
These partial derivatives are also continuous in
$\R\times[\epsilon, T]$ as a result of Remark~\ref{remark:u_c1}.
Then the statement in (ii) follows since the choice of $\epsilon$
is arbitrary.\\
{\bf{(iii)}} It is clear that the functions $x \to J(x,t)$ and
$t\rightarrow J(x,t)$ are nondecreasing functions since $x
\rightarrow u(x,t)+e^x$ and $t\rightarrow u(x,t)$
 are nondecreasing.

  The condition \eqref{eq:condition-nu} means that the support of the measure $\nu$
 is not bounded from above.
  As a result we have that the set $A=\{z : x+z \in \C\}$ has positive measure, i.e., $\nu(A)>0$ for any $x\in \R$.
For any $z \in A$ we have that
   $\partial_x u(x+z,t)+e^{x+z}>0$, which is equivalent to $\partial_S V(Se^z,
 t)+1>0$. The latter follows from the convexity of the function $V$ and \eqref{eq:V_3}.
If $z \notin A$, then clearly  $\partial_x u(x+z,t)+e^{x+z}=0$.
Using these facts in the first equation in \eqref{eq:addt1}, we see that \ref{eq:condition-nu} yields $\partial_x
 J(x,t)>0$ in $\R\times[0,T]$. On the other hand, when $q>0$ the condition assumed on $\nu$ can be dropped.

Moreover, when $x \geq b(t)$ \eqref{eq:condition-q} ensures that
 $\nu(A)>0$. Then \eqref{eq:jt} follows from
 Proposition \ref{prop:u_t positive}.
\end{proof}
In the rest of the paper, we will assume either
\eqref{eq:condition-nu} or $q>0$ and \eqref{eq:condition-q} are
satisfied. Indeed, in the two well-known examples of jump
diffusions, Kou's model and Merton's model (see \cite{cont}
p.111), in which $\nu$ is the double exponential and normal
distribution respectively, condition \eqref{eq:condition-nu} is
fulfilled.

As the consequence of Lemma \ref{lemma:j}, the level curve \begin{equation}\label{eq:B}
 B(t) \triangleq \left\{x : J(x,t) = 0, t\in[0,T]\right\}.
\end{equation}
is well defined.
 $B(0)$, which is the unique
solution of the integral equation,
\begin{equation}\label{eq:B(0)-eq}
 J_0(x) = qe^x - rK +\lambda \int_{\R}\left(e^{x+z}-K\right)^+
 \nu(dz) =0.
\end{equation}
will be crucial in describing the behavior of $b(t)$ close to 0
(see Section~\ref{sec:freebdrybeh}).

\begin{remark}
When $r=0$, Lemma~\ref{lemma:j} (i) implies that
$B(t)=-\infty$. On the other hand, the proof in the following
lemma tell us that $B(t)\geq b(t)$. Therefore $b(t)=-\infty$ in this case. We will assume
$r>0$ in the rest of the paper to exclude this trivial case.
\end{remark}

This level curve $B(t)$ will be crucial in analyzing the
regularity properties of the free boundaries in the rest of this
section. Let us analyze its properties first.

\begin{lemma}\label{lemma:level-B}
 \begin{enumerate}
  \item $B(t)$ is non-increasing,
  \item $B(t) \in C^1(0,T]\cap C[0,T]$,
  \item $B(t) > b(t)$ for $t\in(0,T]$. Here $b(t)$ is the free
  boundary in \eqref{eq:u_1} - \eqref{eq:u_4}.
 \end{enumerate}
\end{lemma}
\begin{proof}
 {\bf{(i)}} The proof follows from Lemma \ref{lemma:j} (iii).
 \\ {\bf{(ii)}} We have the continuity of $B$ because $J(x,t)$ is continuous and strictly increasing in $x$ (see Lemma~\ref{lemma:j} (ii) and (iii)). Let us focus on the differentiability in the following.
 It follows from
 Lemma~\ref{lemma:j} (ii) that $J(x,t)$ is a $C^1$ function in
 $\R\times(0, T]$. Moreover, it follows from \eqref{eq:jt} and $B(t)\geq b(t)$ (which we will prove in the Step 1 in (iii)) that
 \[
  \left.\partial_t J(x,t_0)\right|_{x=B(t_0)} >0, \quad t_0 \in
  (0,T_0].
 \]
 Therefore, the Implicit Function Theorem implies that there
 exists an open set $U$ containing $t_0$ such that
 \[
  B(t) \in C^1(U).
 \]
 Then the statement in (ii) follows after pasting different
 neighborhoods for all points $t\in(0,T]$ together.

 {\bf{(iii)}} The proof consists of two steps:

 {\bf{Step 1:}} First we show that $B(t)\geq b(t)$.
 If these is a $t_0\in(0,T]$ such that $B(t_0) < b(t_0)$, from
 the definition of $B(t)$ and the fact that $x\rightarrow J(x,t)$ is strictly
 increasing, we obtain $J(x,t_0) >0$ for all $x\in(B(t_0),
 b(t_0))$. Combining with \eqref{eq:lu=j}, we have
 \[
  \L u(x,t_0) <0, \quad \text{ for any } x\in (B(t_0), b(t_0)),
 \]
 which contradicts with \eqref{eq:u_5}.

 {\bf{Step 2:}} Second, we show that $B(t) \neq b(t), t\in (0,T]$.
Since $b(t) < \log{K}$ (thanks to Lemma~\ref{lemma:b_strictly}) and
 $t \to B(t)$ is non-increasing, it is clear that $B(t)>b(t)$ for any $t\in (0,
 t^*)$ where $t^* = T \wedge \sup \{t\in \R_+ : B(t)=\log{K}\}$. Hence we only need to focus on the region where $B(t)<\log{K}$. If there is a $t_0\in(0,T]$ such that $B(t_0)=b(t_0)$, we will
 derive a contradiction in the following.

 First, let us define the
 region $\Omega \triangleq \left\{(x,t) \,|\, B(t) < x < \log{K}, t\in
 (0,T]\right\}$. Because of the result in Step 1, $\Omega\subset \C$.
 Hence $u(x,t)$ satisfies
 \[
  \L_{\D} u(x,t) = \lambda \int_{\R}
  u(x+z,t)\nu(dz), \quad (x,t)\in \Omega.
 \]
Let us define $\xi \triangleq x-B(t)$, $\tilde{u}(\xi,t)
 \triangleq
 u(x,t)$ and $\tilde{g}(\xi,t) \triangleq (K-e^{\xi+B(t)})^+ =
 g(x)$. In the region $\widetilde{\Omega} \triangleq \{(\xi,t) |\, 0< \xi < \log{K}- B(t), t\in(0,T]\}$ we have
 \begin{equation}\label{eq:tltu}
  \tilde{\L}_{\D} \tilde{u} \triangleq \frac{\partial
  \tilde{u}}{\partial t} - \frac12 \sigma^2 \frac{\partial^2 \tilde{u}}{\partial
  \xi^2} - \left(\mu + B'(t)
  -\frac12 \sigma^2\right)\frac{\partial \tilde{u}}{\partial \xi} + (r+\lambda) \tilde{u}
  = \lambda \int_{\R} \tilde{u}(\xi+z,t)\nu(dz).
 \end{equation}
 since $B(t)\in C^1(0,T]$.
 On the other hand,
 \begin{equation}\label{eq:tltg}
 \begin{split}
  \tilde{\L}_{\D} \tilde{g} &= - e^{\xi+B(t)}B'(t) + \frac12 \sigma^2 e^{\xi+B(t)} +
  \left(\mu + B'(t) -\frac12 \sigma^2\right)e^{\xi+B(t)} +
  (r+\lambda) \left(K- e^{\xi+B(t)}\right)\\
  & = - \left[qe^{\xi+B(t)} - rK +\lambda \int_{\R} \left(e^{\xi+B(t)+z}-K\right)\nu(dz)\right].
 \end{split}
 \end{equation}
 Therefore, we obtain from
 \eqref{eq:tltu} and \eqref{eq:tltg} that
 \begin{equation}\label{eq:tltu-tg}
  \tilde{\L}_{\D} \left(\tilde{u}-\tilde{g}\right)(\xi,t)  =
  qe^{\xi+B(t)} - rK + \lambda \int_{\R} \left[\tilde{u}(\xi+z,t) + e^{\xi+B(t)+z}
  -K\right]\nu(dz) = J(\xi+B(t),t),
 \end{equation}
for $(\xi,t) \in \widetilde{\Omega}$.
 Note that  $J(x,t) >0$ when $x>B(t)$. Therefore
 \eqref{eq:tltu-tg} yields
 \begin{equation}\label{eq:tl-pos}
  \tilde{\L}_{\D} \left(\tilde{u}-\tilde{g}\right) >0, \quad
  (\xi,t)\in \widetilde{\Omega}.
 \end{equation}

 On the other hand, from our assumption $\xi_0 \triangleq
 b(t_0)-B(t_0)=0$. Moreover, there clearly exists a ball $\mathcal{B}\subset
 \widetilde{\Omega}$ such that 1) $\overline{\mathcal{B}}\cap \{\xi=0\} =
 (\xi_0,t_0)$; 2) $(\tilde{u}-\tilde{g})(\xi,t) > (\tilde{u}-\tilde{g})(\xi_0,t_0) =0$ for all $(\xi,t)\in \mathcal{B}$, since
$(\tilde{u}-\tilde{g})(\xi,t) =
 (u-g)(x,t)>0$ when $x>B(t)\geq b(t)$. Now applying Hopf's Lemma (see Theorem 17 in page 49 of
 \cite{friedman-1}) to $\tilde{u}-\tilde{g}$ in $\mathcal{B}$, we obtain
 \begin{equation}\label{eq:hopf}
  \frac{\partial}{\partial \xi} \left(\tilde{u}
  -\tilde{g}\right)(\xi_0,t_0) > 0,
 \end{equation}
 which contradicts with the smooth fit property at $(\xi_0, t_0)$, i.e., $\partial_{\xi}(\tilde{u}-\tilde{g})(\xi_0,t_0) = \partial_x
 (u-g)(b(t_0),t_0)=0$.
\end{proof}

\begin{remark}\label{remark:interior-ball}
 In the proof of Lemma \ref{lemma:level-B} (iii), the reason we
 work with the domain $\widetilde{\Omega}$ instead of the domain
 $\Omega$ is that $\Omega$ may not satisfy the interior ball
 condition (see Theorem 17 in page 49 of \cite{friedman-1}),
 which is a crucial assumption of the Hopf Lemma. If one can show $B(t)\in
 C^2$, the interior ball condition automatically holds for $\Omega$ (see Remark in page 330 of
 \cite{evans}). However, $B(t)\in C^2$ does not follow directly from the Implicit Function Theorem, because $J(x,t)$ is not expected to be a $C^2$
 function in a neighborhood of the point $(b(t_0),t_0)$, for any $t_0$, as a
 result of the discontinuity of
 $\partial^2_{x}u(x,t)$ across the free
 boundary $b(t)$ (see the following corollary).
\end{remark}

As a corollary of Lemma \ref{lemma:level-B} (iii), $\partial^2_x
u(x,t)$ does not cross the free boundary continuously.
\begin{corollary}\label{cor:u_xx_pos}
\begin{equation}\label{eq:u_xx_pos}
 \frac{\partial^2}{\partial x^2} u\left(b(t)+,t\right) \triangleq \lim_{x\downarrow b(t)} \frac{\partial^2}{\partial x^2} u(x,t) >
 - e^{b(t)}, \quad t\in(0,T].
\end{equation}
(This is equivalent to $\lim_{S\downarrow
s(t)}\partial^2_{S}V(S,t)>0, t\in[0,T)$.)
\end{corollary}
\begin{proof}
 On the one hand, since $B(t)>b(t)$ and $x\rightarrow J(x,t)$ is strictly increasing, we have
 \begin{equation}\label{eq:j-pos}
 J(b(t),t)<0, \quad t\in(0,T],
 \end{equation}
 On the other hand, from the continuity of $u$, \eqref{eq:u_3}, \eqref{eq:u_1} and Proposition~\ref{prop:u_t},
it follows that
 \begin{equation}\label{eq:j-cont}
 \begin{split}
  0 &= \lim_{x\downarrow b(t)} \L u(x,t) = - \frac12 \sigma^2 \lim_{x\downarrow
  b(t)} \frac{\partial^2}{\partial x^2} u(x,t) - \frac12 \sigma^2 e^{b(t)} -
  \left\{qe^{b(t)}-rK + \lambda\int_{\R}\left[u(b(t)+z,t) + e^{b(t)+z}
  -K\right]\nu(dz)\right\}\\
  & = - \frac12 \sigma^2 \lim_{x\downarrow
  b(t)} \frac{\partial^2}{\partial x^2} u(x,t) - \frac12 \sigma^2 e^{b(t)} - J(b(t),t).
 \end{split}
 \end{equation}
 The inequality \eqref{eq:u_xx_pos} now follows from combining
 \eqref{eq:j-pos} and \eqref{eq:j-cont}.
\end{proof}

\subsection{The behavior of the free boundary close to maturity}\label{sec:freebdrybeh}
We are ready to analyze the regularity of the free boundaries. The
continuity of the free boundaries for differential equations with
or without integral terms have been studied intensively, see e.g.
\cite{friedman-3}, \cite{pham}, \cite{yang} and \cite{lamberton}.
For the American option in jump diffusions, \cite{pham} showed the
continuity of the free boundary under the technical condition
\begin{equation}\label{cond:pham}
 r> q+\lambda \int_{\R_+} \left(e^z-1\right) \nu(dz).
\end{equation}
In \cite{yang}, this condition was removed in the proof of the
continuity. Moreover, in their Theorem 5.3, they showed that
\begin{equation}
b(0+) \triangleq \lim_{t\rightarrow 0^+}b(t) = \min\{\log{K},B(0)\}
= \left\{\begin{array}{ll} \log{K}, & r\geq
q+\lambda\int_{\R_+}(e^z-1)\nu(dz)\\
B(0), & r< q+\lambda\int_{\R_+}(e^z-1)\nu(dz)
\end{array}\right., \label{eq:b(0)}
\end{equation}
in which $B(0)$ is the unique solution of \eqref{eq:B(0)-eq}. The
same result has been shown for the exponential L\'{e}vy models in
\cite{lamberton}.

\subsection{H\"{o}lder continuity of the free boundary}
In the following, the function $J_0(x)$ in \eqref{eq:def-J0} and
the Maximum Principle will play a crucial role in showing that $t
\to b(t)$ is H\"{o}lder continuous.

\begin{lemma}\label{lemma:b_holder_est}
Let $b(t)$ be the free boundary in Lemma~\ref{lemma:b_strictly}.
For any $\epsilon>0$, if there exists $\delta>0$ such that for any
$t_1$ and $t_2$ satisfying $\epsilon \leq t_1 < t_2\leq T$ and
$t_2-t_1\leq \delta$  one has
\begin{equation}\label{eq:ineq_u_bt1}
u(b(t_1), t) - u(b(t_1), t_1) \leq C_{\epsilon}(t_2-t_1)^{\alpha},
\quad t_1\leq t\leq t_2,
\end{equation}
in which $0< \alpha \leq 1$ and $C_{\epsilon}$ is a constant
that does not depend on $t_1$ and $t_2$, then there exists $\delta'\in
(0,\delta]$ such that
\begin{equation}\label{eq:b_holder_al}
b(t_1)-b(t_2) \leq C'_{\epsilon} (t_2-t_1)^{\frac{\alpha}{2}},
\quad 0 \leq t_2-t_1\leq \delta',
\end{equation}
in which $C'_{\epsilon}$ is another positive constant that is
independent of $t_1$ and $t_2$.
\end{lemma}

\begin{proof}
This proof is motivated by Lemma 5.1 in \cite{friedman-shen}. For
any $t_1$ and $t_2$ such that $\epsilon \leq t_1<t_2\leq T$ and
$t_2-t_1\leq \delta$, let us consider the domain $D \triangleq
\{(x,t): b(t)< x < b(t_1), t_1<t<t_2\}$. (In what follows, we will
choose $t_1$ and $t_2$ close to each other, i.e. we will find an
appropriate $\delta'$ such that $t_2-t_1 \leq \delta'$.) Let
$\overline{D}$ be the closure of the domain $D$.

In the following, we will show that the function
\begin{equation}\label{eq:def-chi}
\chi(x) =
\left\{\left[\sqrt{C_{\epsilon}}(t_2-t_1)^{\frac{\alpha}{2}} +
\beta(x-b(t_1))\right]^+\right\}^2, \quad b(t_2)\leq x\leq b(t_1)
\end{equation}
satisfies $\chi(x)\geq (u-g)(x,t)$ on the domain $D$ for suitably
chosen positive constant $\beta$.

It is clear that $\chi(x)=0$, when $x\leq b(t_1)-
\frac{\sqrt{C_{\epsilon}}}{\beta}(t_2-t_1)^{\frac{\alpha}{2}}
\triangleq \xi$. We also have $\chi(b(t_1))=C_{\epsilon}
(t_2-t_1)^{\alpha} \geq u(b(t_1),t) - g(b(t_1))$ for $t_1\leq
t\leq t_2$ because of the assumption (\ref{eq:ineq_u_bt1}). On the
other hand, $\chi(b(t))\geq 0 = u(b(t),t)-g(b(t))$. Therefore on
the parabolic boundary of the domain $D$, we have that $\chi\geq
u-g$. We will show that this holds for all $(x,t)\in D$. To this
end, we will compare $\L\chi$ with $\L(u-g)$ using the Maximum
Principle. Note that $\chi$ is carefully chosen so that it has a
continuous first derivative and a bounded second derivative. These
properties of $\chi$ makes the application of the Maximum
Principle for weak solutions (see e.g. Corollary 7.4 in
\cite{lieberman}) possible.

First, for $(x,t)\in D$ let
us estimate the integral term:
\begin{equation}\label{eq:est-integral}
\begin{split}
\lambda \int_{\R} \chi(x+z) \nu(dz) &= \lambda \int_{z\geq \xi-x}
\left\{\sqrt{C_{\epsilon}}(t_2-t_1)^{\frac{\alpha}{2}} +
\beta(x+z-b(t_1))\right\}^2 \nu(dz) \\
&\leq \lambda \int_{z\geq \xi-x}
\left\{\sqrt{C_{\epsilon}}(t_2-t_1)^{\frac{\alpha}{2}} + \beta
z\right\}^2 \nu(dz) \\
&\leq 2\lambda \int_{z\geq \xi-x} \left[C_{\epsilon}
(t_2-t_1)^{\alpha} + \beta^2 z^2\right] \nu(dz)
\\
&\leq 2\lambda \left[C_{\epsilon} (t_2-t_1)^{\alpha} + \beta^2
M\right].
\end{split}
\end{equation}
for a sufficiently large positive constant $M$ independent of
$t_1$ and $t_2$. To obtain the first inequality, we used $x<
b(t_1)$ for $(x,t)\in D$. The third inequality follows, because
$\int_{\R} e^z \nu(dz) < +\infty$ in \eqref{eq:xi} and $z$ is
bounded from below.

With the estimate (\ref{eq:est-integral}), we can calculate
$\L\chi$ inside the domain $D$.
\begin{equation}\label{eq:est-LX}
\begin{split}
\L\chi(x) &= \left[-\sigma^2 \beta^2 - \left(\mu -\frac12
\sigma^2\right)2\beta\chi^{\frac12} + (r+\lambda)\chi\right]
1_{\{x>\xi\}} - \lambda \int_{\R} \chi(x+z, t)\nu(dz)
\\
&\geq - \left[\frac{(\mu-
\sigma^2/2)^2}{r+\lambda}+\sigma^2\right]\beta^2 1_{\{x>\xi\}} -
2\lambda \left[C_{\epsilon} (t_2-t_1)^{\alpha} +
\beta^2 M\right] \\
&\geq -E\beta^2 - F(t_2-t_1)^\alpha,
\end{split}
\end{equation}
in which $E \triangleq \frac{(\mu
-\sigma^2/2)^2}{r+\lambda}+\sigma^2 + 2\lambda M$ and $F
\triangleq 2\lambda C_{\epsilon}$ are positive constants.

Recall that for any $\eps>0$, $b(\epsilon) < \min\{\log{K}, B(0)\}$
and that the strictly increasing function $J_0$ defined in
(\ref{eq:def-J0}) satisfies $J_0(x)<0$ for $x<B(0)$. Using these
observations and (\ref{eq:lg}) it can be seen that for any $x\leq
b(\epsilon)$ we have
\begin{equation}\label{eq:est-Lg}
 \L g(x) = -J_0(x) \geq - J_0(b(\epsilon)) >0.
\end{equation}

Now choosing
\begin{equation}\label{eq:chs-c}
c=- J_0(b(\epsilon)) >0
\end{equation}
and $\delta' =
\min\{ \left(\frac{c}{2F}\right)^{1/\alpha}, \delta \}$ and $\beta
\leq \sqrt{\frac{c}{2E}}$, we have that
\[
 \L \chi(x) (x) \geq -c \geq \L(u-g)(x,t), \quad (x,t)\in D.
\]
Considering $\Psi = \chi-u+g$, we have $\L\Psi \geq 0$ in $D$ and
$\Psi\geq 0$ on the parabolic boundary of $D$. It follows from the
Maximum Principle for weak solutions that $\Psi\geq 0$ in $D$, i.e.,
\begin{equation}\label{eq:X>u-g}
 \chi(x)\geq (u-g)(x,t), \quad (x,t)\in D.
\end{equation}
Observe that $(u-g)(x,t)=0$ if $x\leq \xi$. For any $(x,t)\in D$,
since $(u-g)(x,t)>0$, we can see that $x>\xi$. This gives us
\begin{equation}\label{eq:b-holder-al-1}
\inf_{t_1\leq t \leq t_2} b(t) \geq b(t_1)
-\frac{\sqrt{C_{\epsilon}}}{\beta} (t_2-t_1)^{\frac{\alpha}{2}},
\quad 0 < t_2-t_1\leq \delta'.
\end{equation}
We have shown the free boundary $b(t)$ is continuous and strictly
decreasing in Lemma~\ref{lemma:b_strictly}. Along with this fact,
the inequality (\ref{eq:b-holder-al-1}) gives us
(\ref{eq:b_holder_al}) with $C'_{\epsilon} =
\sqrt{C_{\epsilon}}/\beta$.
\end{proof}

Now we are ready to state the main result of this section.
\begin{theorem}\label{theorem:holder-cont}
 Let $b(t)$  be the free boundary in problem (\ref{eq:u_1}) - (\ref{eq:u_4}), then for any $\epsilon>0$ if $\epsilon \leq t_1 < t_2 \leq
 T$, and $t_2-t_1$  is sufficiently small, then
\begin{equation}
 b(t_1)-b(t_2) \leq C_{\epsilon}(t_2-t_1)^{\frac58}, \quad ,
\end{equation}
in which $C_{\epsilon}$ is a positive constant independent of
$t_1$ and $t_2$.
\end{theorem}

\begin{proof}
The proof will follow by applying Lemma~\ref{lemma:b_holder_est}
twice. The first application will show that
 $b(t)$ is H\"{o}lder continuous with exponent
$\frac12$. Applying   Lemma~\ref{lemma:b_holder_est}  for the second time we will upgrade the H\"{o}lder exponent to
$\frac58$.

As a result of Propositions~\ref{prop:u_t} and \ref{prop:u_t
positive} for any $\epsilon>0$, $t_1$ and $t_2$ satisfying
$\epsilon \leq t_1<t_2 \leq T$ we have that
\begin{equation} \label{eq:est-up-1}
u(b(t_1),t) - u(b(t_1), t_1) \leq \max_{t_1\leq s\leq t}
\frac{\partial u}{\partial t}(b(t_1),s) (t-t_1) \leq C_1
(t_2-t_1),
\end{equation}
where $C_1 = \max_{\epsilon\leq s\leq T} \partial_t u(b(t_1),s)$
is a positive constant. Now as a result of
Lemma~\ref{lemma:b_holder_est}, we know that there exists a
sufficiently small constant $\delta_1 \in (0,T-\epsilon]$ such
that
\begin{equation}\label{eq:b-holder-1/2}
 b(t_1) - b(t_2) \leq C'_1(t_2-t_1)^{\frac12}, \quad 0\leq t_2-t_1
 \leq \delta_1,
\end{equation}
in which $C'_1$ is a positive constant that does not depend on
$t_1$, $t_2$ and $\delta_1$.

It follows from Lemmas 2.8 and 2.11 in \cite{yang} and the Sobolev
Embedding Theorem (see also (\ref{eq:hold-unt}) in Appendix A.3)
that for any $a<b<\log{K}$ and $t \in [t_1,t_2]$,
\begin{equation}\label{eq:ut-holder}
\left|\frac{\partial u}{\partial t}(x,t) - \frac{\partial
u}{\partial t} (\overline{x},t)\right| \leq \tilde{C}
\left|x-\overline{x}\right|^{\frac12}, \quad x,\, \overline{x} \in
(a,b),
\end{equation}
in which $\tilde{C}$ is a positive constant that does not depend
on $t$. Taking $x=b(t_1)$ and $\overline{x}=b(t)$ in
(\ref{eq:ut-holder}) and using Proposition~\ref{prop:u_t}, we
obtain
\begin{equation}\label{eq:est-up-2}
 0\leq \frac{\partial u}{\partial t}(b(t_1),t) \leq
 \tilde{C}\left|b(t_1)-b(t)\right|^{\frac12} \leq \tilde{C}
 |b(t_1)-b(t_2)|^{\frac12}, \quad t_1\leq t\leq t_2,
\end{equation}
where the third inequality follows from $b(t)$ being strictly
decreasing in Lemma~\ref{lemma:b_strictly}. Combining
(\ref{eq:b-holder-1/2}) and (\ref{eq:est-up-2}), we get
\begin{equation}\label{eq:est-up-3}
0\leq \frac{\partial u}{\partial t}(b(t_1), t) \leq C_2
(t_2-t_1)^{\frac14}, \quad t_1\leq t\leq t_2, \, 0\leq t_2-t_1\leq
\delta_1.
\end{equation}
As a result
\begin{equation}
 u(b(t_1),t)- u(b(t_1),t_1) \leq \max_{t_1\leq s\leq t_2} \frac{\partial
 u}{\partial t}(b(t_1),s)(t_2-t_1) \leq
 C_2(t_2-t_1)^{\frac54}.
\end{equation}
Applying Lemma~\ref{lemma:b_holder_est} for the second time, we know that there
exists $\delta_2 \in (0, \delta_1]$ such that
\begin{equation}\label{eq:b-holder-5/8}
 b(t_1)-b(t_2) \leq C_{\epsilon}(t_2-t_1)^{\frac58}, \quad 0\leq t_2-t_1\leq
 \delta_2,
\end{equation}
where $C_{\epsilon}$ is a positive constant that does not depend on $t_1$, $t_2$
and $\delta_2$.
\end{proof}

\section{The free boundary is continuously differentiable}\label{sec:cont-diff}
In this section, we will investigate the continuous
differentiability of the free boundary. In Theorem 5.6 in
\cite{yang}, the authors have shown that $b(t)\in C^1(0,T]$, with the
extra condition
\begin{equation}\label{eq:condition}
 r \geq q+\lambda \int_{\R_+} \left(e^z-1\right) \nu(dz).
\end{equation}
Thanks to Corollary~\ref{cor:u_xx_pos} and
Theorem~\ref{theorem:holder-cont}, we can show the continuous
differentiability of the free boundary without imposing this extra
condition.
\begin{remark}
 If condition (\ref{eq:condition}) is not satisfied, we can
 see from (\ref{eq:b(0)}) that there is a gap between $\lim_{t\rightarrow 0^+}
 b(t)$ and $b(0)=\log{K}$. Therefore it is impossible to have $b(t)$ to be even continuous at $t=0$. But we shall see that it is continuously differentiable for all $t \in (0,T]$.
\end{remark}

Let us consider the time derivative $\partial_tu(x,t)$. Recall
that $u(x,t)$ is the solution of (\ref{eq:u_1}) - (\ref{eq:u_4}).
Using the assumption \eqref{eq:ass-sigma}, the time derivative
$w=\partial_t u(x,t)$ satisfies the following partial differential
equation
\begin{equation}\label{eq:w_1}
\begin{split}
 &\L_{\D} w= h(x,t), \quad x>b(t), \; t\in (0,T],\\
 &w(b(t),t) =0, \quad \lim_{x\rightarrow +\infty} w(x,t)=0, \quad t\in(0,T], \\
 &w(x,0)=\lim_{t\rightarrow 0}\partial_t u(x,t), \quad x\geq b(0),
\end{split}
\end{equation}
in which
\begin{equation}\label{eq:def-h}
h(x,t) \triangleq \lambda\int_{\R} \frac{\partial}{\partial
t}u(x+z,t) \nu(dz) + \sigma \cdot (\partial_t \sigma) \cdot
\left(\frac{\partial^2 u}{\partial x^2} - \frac{\partial
u}{\partial x}\right).
\end{equation}
When $x<b(t)$, we also have $w(x,t)=0$. Given $u(x,t)$ and $b(t)$,
(\ref{eq:w_1}) is a parabolic differential equation for $w(x,t)$.
In this equation, the boundary conditions for $w(x,t)$ along
$b(t)$ and at the infinity follow from Proposition \ref{prop:u_t}
and Lemma \ref{lemma:u_t_inf}.

In order to show the differentiability of the free boundary, we
need to study the behavior of $\frac{\partial^2}{\partial
x\partial t}u$ at the free boundary (by first making sure that the
cross derivatives exist in the classical sense), which is carried
out in the following lemma.

\begin{lemma}\label{lemma:u_xt cont}
\begin{itemize}
\item[(i)] As a function of $t$, $\frac{\partial^2}{\partial
x\partial t}u(b(t)+,t) \triangleq \lim_{x\downarrow
b(t)}\frac{\partial^2}{\partial x\partial t} u(x,t)$ is continuous
on $(0,T]$. \item[(ii)] Moreover, the function
$\frac{\partial^2}{\partial x\partial t}u(x,t)$ is continuous for
$x>b(t)$, $t\in (0,T]$.

\end{itemize}
\end{lemma}
This lemma is a slight generalization of the result in
\cite{Cannon} to the parabolic integro-differential equation
\eqref{eq:w_1}. Considering the integral term $h$ in
(\ref{eq:w_1}) as the driving term, this lemma follows from using
the same technique presented in Section 1 of Chapter 8 in
\cite{friedman-1}. We will postpone this proof to the Appendix
\ref{app:A.2}. We are now ready to state and prove the main
theorem of this section.
\begin{theorem}\label{theorem:cont-diff}
Let $b(t)$ be the free boundary in the boundary value problem
(\ref{eq:u_1}) - (\ref{eq:u_4}), then $b(t)\in C^1(0,T]$.
\end{theorem}
\begin{proof}
First, we will show $b(t)$ is differentiable at $t_0\in (0,T]$.
Let us define $\rho=\partial^2_x u(b(t_0)+,t_0) + e^{b(t_0)}$.
Corollary~\ref{cor:u_xx_pos} implies that $\rho>0$.

For sufficiently small $\epsilon>0$, it follows from
(\ref{eq:u_3}) that
\[
 \frac{1}{\epsilon} \left[\frac{\partial}{\partial x}u(b(t_0), t_0) - \frac{\partial}{\partial x}u(b(t_0-\epsilon), t_0-\epsilon) + e^{b(t_0)} -
 e^{b(t_0-\epsilon)}\right]=0.
\]
Applying the Mean Value Theorem yields
\begin{equation}\label{eq:lim-u_t}
\left(\frac{\partial^2}{\partial x^2}u(b(t_0)+y, t_0)+
e^{b(t_0)+y}\right)\frac{b(t_0)-b(t_0-\epsilon)}{\epsilon} =
-\frac{\partial^2}{\partial x\partial t}
u(b(t_0-\epsilon),t_0-\tau),
\end{equation}
for some $ y \in (0,b(t_0-\epsilon)-b(t_0))$ and $\tau \in (0,\epsilon)$.
Letting $\epsilon \rightarrow 0$ in
(\ref{eq:lim-u_t}) and using Lemma \ref{lemma:u_xt cont} (ii), we
obtain
\begin{equation}\label{eq:b'(t)}
\lim_{\epsilon\rightarrow 0}
\frac{b(t_0)-b(t_0-\epsilon)}{\epsilon} =
-\frac{\frac{\partial^2}{\partial x\partial
t}u(b(t_0)+,t_0)}{\frac{\partial^2}{\partial
x^2}u(b(t_0)+,t_0)+e^{b(t_0)}},
\end{equation}
which implies that $b(t)$ is differentiable since $\rho>0$.
Moreover, from (\ref{eq:u_1_para}) and Proposition \ref{prop:u_t},
we have
\[
 \frac{\partial^2}{\partial x^2}u(b(t)+,t)=
 \frac{2(r+\lambda)}{\sigma(b(t),t)^2}K +
 \left(\frac{2(\mu-r-\lambda)}{\sigma(b(t),t)^2}-1\right)e^{b(t)} - \frac{2}{\sigma(b(t),t)^2} f(b(t),t),
\]
which is clearly a continuous function of $t$ on $t\in (0,T]$,
since $b(t)$ is a continuous function and $\sigma(x,t)$ is
continuous from our assumption \eqref{eq:ass-sigma}. Along with
Lemma \ref{lemma:u_xt cont} (i), we can see from (\ref{eq:b'(t)})
that $b(t)\in C^1(0,T]$.
\end{proof}

\section{Higher order regularity of the free boundary} \label{sec-inf-dif}
In the previous section, we have proved that the free boundary
$b(t)$ is continuously differentiable. In this section, we will
upgrade their regularity. Throughout this section, for the
simplicity of the notation, we will assume that $\sigma$ is a
positive constant. In this case, $h(x,t)=\lambda\int_{\R}
\frac{\partial}{\partial t}u(x+z,t) \nu(dz)$, which is bounded
thanks to Lemma~\ref{lemma:ut-bounded}. More generally, if
$\sigma=\sigma(x,t)$, $h(x,t)$ is given in \eqref{eq:def-h}. If we
assume $\sigma(x,t)\in C^{\infty}(\R\times [0,T])$ with all its
derivatives bounded and $\delta \leq \sigma \leq \Delta$ for some
positive constants $\delta$ and $\Delta$, the same arguments in
this section can still be carried through. Because of
Lemmas~\ref{lemma:lip-holder} and \ref{lemma:ut-bounded}, we can
see from the equation \eqref{eq:u_1} that $\partial^2_x u(x,t)$ is
also bounded in $\R\times[\epsilon, T]$ for any $\epsilon>0$.
Hence, $h(x,t)$ is also bounded in this general case.

First, let us derive an identity for $b'(t)$.  Since $b(t)$ is
differentiable, taking derivative with respect to $t$ on both
sides of (\ref{eq:u_3}), we have
\begin{equation}\label{eq:u_3_t}
\frac{\partial^2}{\partial x^2}u(b(t)+,t) b'(t) +
\frac{\partial^2}{\partial x\partial t}u(b(t)+,t) =
-e^{b(t)}b'(t).
\end{equation}
The term $\partial_x^2 u(b(t)+,t)$ can be represented as
\begin{equation}\label{eq:u_xx}
 \frac{\partial^2}{\partial x^2}u(b(t)+,t) =
 \left(\frac{2(\mu-r-\lambda)}{\sigma^2}-1\right)e^{b(t)} +
 \frac{2(r+\lambda)}{\sigma^2}K - \frac{2}{\sigma^2} f(b(t),t).
\end{equation}
Plugging (\ref{eq:u_xx}) back into (\ref{eq:u_3_t}) and recalling
$w=\partial_t u$, we obtain
\begin{equation}\label{eq:def-b'}
 b'(t)= - \frac{\frac{\sigma^2}{2}\frac{\partial}{\partial
 x}w(b(t)+,t)}{(\mu-r-\lambda)e^{b(t)}+ (r+\lambda)K - f(b(t),t)}, \quad t\in (0,T].
\end{equation}

We can see from equations (\ref{eq:w_1}) that $w(x,t)$ is the
solution of a formal Stefan problem in the unbounded continuation
regions $\mathcal{C}$. \cite{schaeffer} gave a proof of the
infinite differentiability of the free boundary of a one
dimensional Stefan problem in a bounded domain. By introducing the
new variable $\xi=\frac{x}{b(t)}$, he reduced the  problem into a
fixed boundary problem on a bounded domain. However, if we apply
the same change of variables we will have unbounded coefficients
in the corresponding fixed boundary problem. Instead, similar to
the change of variables in the proof of Lemma~\ref{lemma:level-B}
(iii), we will define
\[
\xi\triangleq x-b(t), \quad v(\xi,t)\triangleq w(x,t),
\]
in which $b(t)$ is the free boundary in (\ref{eq:u_1}) -
(\ref{eq:u_4}). The function $v(\xi,t)$ satisfies the following
fixed boundary equation,
\begin{eqnarray}
&& \frac{\partial v}{\partial t} - \frac12 \sigma^2
\frac{\partial^2 v}{\partial \xi^2} - \left(\mu+ b'(t) -\frac12
\sigma^2 \right)\frac{\partial v}{\partial \xi} + (r+\lambda)v =
h(\xi+b(t),t),\quad (\xi,t)\in (0,+\infty)\times(0,T],
\label{eq:v_fix_1}\\
&& v(0,t)= 0, \quad t\in(0,T],\label{eq:v_fix_2} \\
&& v(\xi,0) = w(\xi+b(0),0), \quad \xi\geq 0. \label{eq:v_fix_3}
\end{eqnarray}
Moreover, we have the following identity
\begin{equation}
b'(t)= - \frac{\frac{\sigma^2}{2}\frac{\partial}{\partial
 \xi}v(0,t)}{(\mu-r-\lambda)e^{b(t)}+ (r+\lambda)K -
 f(b(t),t)}, \quad t\in (0,T].\label{eq:v_fix_4}
\end{equation}

\begin{remark}\label{remark:uniqueness}
Since $b(t)\in C^1(0,T]$, so for any $\epsilon>0$, $b'(t)$ is
continuous and bounded in $[\epsilon,T]$. On the other hand, since
$\partial_t u$ is bounded by Lemma \ref{lemma:ut-bounded}, so
$h(\xi+b(t),t)=\lambda\int_{\R} \partial_t u(\xi+b(t)+z,t)\nu(dz)$
is also bounded when $(\xi,t)\in [0,+\infty)\times[\epsilon,T]$.
As a result, it follows from Theorem 2.6 in page 19 of \cite{lad}
that the parabolic differential equation (\ref{eq:v_fix_1}) with
the initial condition $v(\xi, \epsilon)= w(\xi+b(\epsilon),
\epsilon)$ instead of \eqref{eq:v_fix_3} has at most one bounded
classical solution. It follows from the proof of
Lemma~\ref{lemma:u_xt cont} (i) that $\partial_t u(x,t)$ is a
bounded classical solution, so it is the unique bounded solution
of (\ref{eq:v_fix_1}).
\end{remark}

The following result for parabolic differential equations will be
an essential tool in the proof of the main result in this section.
\begin{lemma}\label{lemma:regular-improve}
Let us assume $w(\xi,t)\in H^{2\alpha,\alpha}([0,+\infty)\times
[\delta,T])$ (for some $\alpha$ and $\delta>0$) satisfies the
following equation
\begin{eqnarray}
 &&\frac{\partial w}{\partial t} - a \frac{\partial^2 w}{\partial
 \xi^2} + \ell\frac{\partial w}{\partial \xi} + cw = d\int_{\R}
 \phi(\xi+z,t) \nu(dz), \quad (\xi,t)\in ((0,+\infty)\times (\delta,
 T))\label{eq:upgrade-v-1}\\
 && w(0,t)=g(t), \; t\in [\delta,T]. \label{eq:upgrade-v-2}
\end{eqnarray}
We assume that  $d \int_{\R}\phi (\xi+z,t) \nu(dz)\in
H^{2\alpha,\alpha}([0,+\infty)\times[\delta,T])$ and that
coefficients $a, \ell, c$ also belong to
$H^{2\alpha,\alpha}([0,+\infty)\times [\delta,T])$ with $\delta
\leq a \leq \Delta$ for some positive constants $\delta$ and
$\Delta$, moreover $g(t)\in H^{1+\alpha}([\delta,T])$. Then
 $w(\xi,t)\in
H^{2+2\alpha,1+\alpha}([0,+\infty)\times[\delta',T]$, for any
$\delta'>\delta$.
\end{lemma}
\begin{proof}
Consider a cut-off function $\eta(t)\in C_0^{\infty}((0,T])$, such
that $\eta(t)=0$ when $t\in (0,\delta]$ and $\eta(t)=1$ for $t\in
[\delta',T]$. The function $\widetilde{w}(\xi,t)=\eta(t)w(\xi,t)$
satisfies
\begin{eqnarray*}
 &&\frac{\partial \widetilde{w}}{\partial t}- a \frac{\partial^2 \widetilde{w}}{\partial
 \xi^2} + \ell \frac{\partial \widetilde{w}}{\partial \xi} +
 c\widetilde{w}=d\int_{\R}\eta(t)\phi(\xi+z,t)\nu(dz) + \frac{\partial \eta}{\partial
 t}w(\xi,t), \quad (\xi,t)\in (0,+\infty)\times(\delta,T],\\
 && \widetilde{w}(0,t)=\eta(t)g(t), \quad t\in [\delta,T],\\
 && \widetilde{w}(\xi,\delta)=0, \quad \xi\geq 0.
\end{eqnarray*}
From our assumptions we have that
\begin{eqnarray*}
 &&d\int_{\R} \eta(t)\phi(\xi+z,t) \nu(dz) + \frac{\partial \eta}{\partial
 t}w(\xi,t) \in H^{2\alpha, \alpha}([0,+\infty)\times[\delta,T]),\\
 &&\eta(t)g(t)\in H^{1+\alpha}([\delta,T]).
\end{eqnarray*}
Moreover, the coefficients of the above differential equation are
all inside space
$H^{2\alpha,\alpha}([0,+\infty)\times[\delta,T])$. In addition,
this equation is uniformly parabolic as the result of $0<\delta
\leq a\leq \Delta$. It follows from regularity estimation for
parabolic differential equations (see Theorem 5.2 in page 320 of
\cite{lad}) that $\widetilde{w}(\xi,t)\in
H^{2+2\alpha,1+\alpha}([0,+\infty)\times [\delta,T])$, which
implies $w(\xi,t)\in H^{2+2\alpha, 1+\alpha}([0,+\infty)\times
[\delta',T])$ by the choice of $\eta(t)$.
\end{proof}

\begin{remark}\label{remark:non-smooth initial cond}
We will apply the previous lemma to $w(x,t)=\partial_t u(x,t)$.
Because the initial condition for $w(x,t)$,
 $\lim_{t\rightarrow 0} \partial_t u(x,t)$, is not smooth, we can not apply Theorem 5.2 in page 320 of \cite{lad} to upgrade the regularity of $w$ directly. This
is the reason we work with $\widetilde{w}$ in the proof of the
previous lemma.
\end{remark}

In order to apply Lemma~\ref{lemma:regular-improve} to
(\ref{eq:v_fix_1}) - (\ref{eq:v_fix_4}), we need H\"{o}lder
continuous coefficients and value functions. Let us first show
that the coefficients in equation (\ref{eq:v_fix_1}) are
H\"{o}lder continuous.

\begin{lemma}\label{lemma:b'-holder}
Let $b(t)$ be the free boundary in (\ref{eq:u_1}) -
(\ref{eq:u_4}). Then $b(t)\in H^{1+\alpha}([\delta,T])$ with
$0<\alpha<\frac12$ for any $\delta>0$.
\end{lemma}
\begin{proof}
 For any $\delta>0$, since $b(t)\in C^1(0,T]$ by Theorem \ref{theorem:cont-diff}, the coefficients
 in equation (\ref{eq:v_fix_1}) are bounded and continuous in
 $[\delta,T]$. On the other hand, because $\partial_t u(x,t)$
 is bounded in $\R\times[\delta,T]$ by Lemma \ref{lemma:ut-bounded}, the function $h(\xi+b(t),t)= \lambda\int_{\R} \frac{\partial}{\partial t}
 u(\xi+b(t)+z,t)\nu(dz)$ is also bounded when $(\xi,t)\in
 [0,+\infty)\times[\delta,T]$. It follows from
 Theorem 9.1 in page 341 of \cite{lad} that equation (\ref{eq:v_fix_1}) has a
 unique solution $v(\xi,t)\in W_q^{2,1}([0,M]\times
 [\delta,T])$ for any $q>1$ and $M>0$.

 By the Sobolev Embedding Theorem (see, for example, Theorem 2.1 in page 61 of \cite{lad}), for
 $q>3$, we have $v(\xi,t)\in H^{\beta,
 \beta/2}([0,M]\times[\delta,T])$ with $\beta= 2-\frac3q$ ($1<\beta<2$). as a result, we have
 \begin{equation}\label{eq:reg-v0}
 \frac{\partial}{\partial \xi} v(0,t) \in
 H^{\frac{\beta-1}{2}}([\delta,T]), \quad \text{with } 0 <
 \frac{\beta-1}{2}<\frac12.
 \end{equation}
Let us analyze the terms in the denominator on the right hand side
of
 (\ref{eq:v_fix_4}). We have that $b(t)\in C^1([\delta,T])$ and that
 \[
   \quad f(b(t),t) =
  \lambda\int_{\R} u(b(t)+z,t) \nu(dz) \in
  C^1([\delta,T]),
 \]
since $u(x,t)\in C^1(\R\times[\delta, T])$ (see
Remark~\ref{remark:u_c1}). Moreover, this
 denominator is also bounded away from 0, because
 \[
  (\mu-r-\lambda){\sigma^2}e^{b(t)} +
  (r+\lambda)K - f(b(t),t) = \frac{\sigma^2}{2}\left(\frac{\partial^2}{\partial
  x^2}u(b(t),t) + e^{b(t)}\right)>0, \quad t\in[\delta,T],
 \]
where the last inequality follows from
Corollary~\ref{cor:u_xx_pos}.
 It is clear from (\ref{eq:v_fix_4}) and (\ref{eq:reg-v0}) that, 
 \[
 b'(t)\in H^{\frac{\beta-1}{2}}([\delta,T]).
 \]
\end{proof}

As a corollary of Lemmas~\ref{lemma:regular-improve} and
\ref{lemma:b'-holder}, we can improve the regularity of the
functions $u(x,t)$.
\begin{corollary}\label{cor:ut-holder}
Let $u(x,t)$ be the classical solution of the boundary value
problem (\ref{eq:u_1}) - (\ref{eq:u_4}). Then $u(\xi+b(t),t)\in
H^{2+2\alpha,1+\alpha}([0,+\infty)\times[\delta',T])$ for any
$\delta'>0$, with $\alpha \in (0,1/2)$. 
\end{corollary}
\begin{proof}
Let $\xi=x-b(t)$, $\kappa(\xi,t)=u(x,t)$ and
$\phi(\xi+z,t)=u(\xi+b(t)+z,t)$. Then $\kappa(\xi,t)$ satisfies a
differential equation of the form (\ref{eq:upgrade-v-1}) and
(\ref{eq:upgrade-v-2}) in Lemma \ref{lemma:regular-improve} with
$g(t)=K-e^{b(t)}$ (in fact $\kappa$ satisfies (\ref{eq:v_fix_1})
when $h$ in the driving term is replaced by $f$). Moreover, by
Lemma \ref{lemma:b'-holder}, the coefficients in this equation
(\ref{eq:upgrade-v-1}) are inside space $H^{\alpha}([\delta,T])$
for any $\delta>0$, and $g(t)\in H^{1+\alpha}([\delta, T])$. In
addition, thanks to the assumption \eqref{eq:ass-sigma}, the
equation (\ref{eq:upgrade-v-1}) is uniformly parabolic.

On the other hand, since $u(x,t)$ is uniformly Lipschitz in $x\in
\R$ and uniformly semi-H\"{o}lder continuous in $t\in[0,T]$ (see
Lemma~\ref {lemma:lip-holder}), and $b(t)$ is continuously
differentiable, it is not hard to see that $\int_{\R}
u(\xi+b(t)+z,t)\nu(dz)\in
H^{2\alpha,\alpha}([0,+\infty)\times[\delta,T])$. Moreover,
$u(\xi+b(t),t)\in H^{2\alpha, \alpha}([0,+\infty)\times [\delta,
T])$ again because of Lemma \ref{lemma:lip-holder}. Now, the
statement follows directly from Lemma \ref{lemma:regular-improve}.
\end{proof}

Armed with Lemmas \ref{lemma:regular-improve},
\ref{lemma:b'-holder} and Corollary \ref{cor:ut-holder}, we can
state and prove the main theorem of this section.

\begin{theorem}\label{theorem:higher-reg-b}
Let $b(t)$ be the free boundary in (\ref{eq:u_1}) -
(\ref{eq:u_4}). Assume that $\nu$ has a density, i.e. $\nu(dz)=
\rho(z) dz$. Let $\alpha \in (0,1/2)$. If $\rho(z)$ satisfies
$\int_{-\infty}^u \rho(z) dz \in H^{2\alpha} (\R_-)$, then $b(t)
\in H^{\frac32 + \alpha}([\epsilon, T])$. On the other hand, if
$\rho(z)\in H^{\ell-1+2\alpha} (\R_-)$ for $\ell\geq 1$, then
$b(t) \in H^{\frac32 + \frac{\ell}{2} + \alpha}([\epsilon, T])$,
for any $\epsilon>0$.
\end{theorem}

\begin{proof}
The proof consists of four steps.

{\bf Step 1.}
 From Lemma \ref{lemma:b'-holder} and Corollary
\ref{cor:ut-holder}, we have that $b(t)\in
H^{1+\alpha}([\delta,T])$ and that $u(\xi+b(t),t)\in
H^{2+2\alpha,1+\alpha}([0,+\infty)\times[\delta',T])$ for any
$\delta'>\delta>0$ with $\alpha \in (0,1/2)$,  which implies that
$\partial_t u(\xi+b(t),t)\in H^{2\alpha,\alpha}([0,+\infty)\times
[\delta',T])$ (see Definition \ref{def:holder space}).

{\bf Step 2.} Assume that there is a positive nonintegral real
number $\beta$ with $2\beta  \leq 2\alpha +\ell$, such that
\begin{eqnarray}
&&b(t)\in H^{1+\beta}([\delta, T]), \label{eq:ind-assumption_1}\\
&&\frac{\partial}{\partial t} u(\xi+b(t),t)\in
H^{2\beta,\beta}([0,+\infty)\times [\delta',T]),
\label{eq:ind-assumption_2}\\
&&u(\xi+b(t),t)\in
H^{2+2\beta,1+\beta}([0,+\infty)\times[\delta',T]),
\label{eq:ind-assumption_3}
\end{eqnarray}
for $\delta'>\delta>0$. We will upgrade the regularity exponent
from $\beta$ to $1/2+\beta$, in steps 2 and 3.

Let us analyze $\partial_t u(\xi+ b(t),t)$. For any integers
$r,s\geq 0$, $2r+s <2\beta$, since $\partial_t u(\xi+b(t)+z,t)=0$
when $z\leq-\xi$, we have
\begin{equation}\label{eq:higher-der-ut}
\begin{split}
&\frac{\partial^s}{\partial \xi^s} \frac{\partial^r}{\partial t^r}
\int_{\R} \frac{\partial}{\partial t} u(\xi+b(t)+z,t) \nu(dz) =
\frac{\partial^s}{\partial \xi^s} \frac{\partial^r}{\partial t^r}
\int_{-\xi}^{+\infty} \frac{\partial}{\partial t}
u(\xi+b(t)+z,t)\rho(z)dz \\
&= 1_{\{s\geq 1\}}\sum_{i=0}^{s-1} \frac{\partial^i}{\partial
\xi^i}\frac{\partial^r}{\partial t^r} \left.
\frac{\partial}{\partial t}u(\xi+b(t)+z,t)   \right|_{z\downarrow
-\xi }  \frac{d^{s-1-i}}{d
\xi^{s-1-i}}\rho(-\xi) \\
&\quad + \int_{-\xi}^{+\infty} \frac{\partial^s}{\partial \xi^s}
\frac{\partial^r}{\partial t^r}\frac{\partial }{\partial t}
u(\xi+b(t)+z,t)\rho(z) dz,
\end{split}
\end{equation}
for any $\xi\geq 0$.

When $t$ is fixed, in the following, we will show
\begin{equation}
\frac{\partial^s}{\partial \xi^s}\frac{\partial^r}{\partial t^r}
\int_{\R} \frac{\partial}{\partial t} u(\xi+b(t)+z,t)\nu(dz) \in
H^{2\beta -[2\beta]}([0,+\infty)), \quad \text{for }
2r+s=[2\beta].\label{eq:holder-ut-xi}
\end{equation}
For  any $\xi_1>\xi_2\geq 0$ such that $\xi_1-\xi_2\leq \rho_0$,
we have
\begin{equation}\label{eq:est-int-ut}
\begin{split}
&\left|\frac{\partial^s}{\partial \xi^s} \frac{\partial^r}
{\partial t^r}\int_{\R} \frac{\partial}{\partial t} u(\xi_1+
b(t)+z, t)\nu(dz) - \frac{\partial^s}{\partial \xi^s}
\frac{\partial^r}{\partial t^r}\int_{\R} \frac{\partial}{\partial
t} u(\xi_2+ b(t)+z, t)\nu(dz) \right|\\
&\leq 1_{\{s\geq 1\}} \sum_{i=0}^{s-1}
\left|\frac{\partial^i}{\partial \xi^i}\frac{\partial^r}{\partial
t^r} \left.  \frac{\partial}{\partial t}
u(\xi + b(t)+z,t) \right|_{z\downarrow -\xi }\right|\left| \left. \left. \frac{ d^{s-1-i}}{d\xi^{s-1-i}} \right(\rho(-\xi_1)-\rho(-\xi_2)\right) \right|\\
&\quad + \int_{-\xi_2}^{+\infty} \left| \left. \left.
\frac{\partial^s}{\partial \xi^s} \frac{\partial^r}{\partial
t^r}\frac{\partial}{\partial t}
\right(u(\xi_1+b(t)+z,t)-u(\xi_2+b(t)+z,t)\right)
\right|\rho(z)dz\\
&\quad  + \int_{-\xi_1}^{-\xi_2} \left|\frac{\partial^s}{\partial
\xi^s} \frac{\partial^r}{\partial t^r}\frac{\partial}{\partial t}
u(\xi_1+b(t)+z,t)\right|\rho(z) dz.
\end{split}
\end{equation}
Let us analyze the right hand side of (\ref{eq:est-int-ut}) term
by term. When $s>1$, since $s-1< 2\beta-1\leq 2\alpha+\ell-1$, we have
$\rho(z)\in H^{2\beta-1}(\R_-)$, which implies
\begin{equation}\label{eq:est-int-ut-1}
\begin{split}
&1_{\{s\geq 1\}} \sum_{i=0}^{s-1} \left|
\frac{\partial^i}{\partial \xi^i}\frac{\partial^r} {\partial t^r}
\left. \frac{\partial}{\partial t} u(\xi+ b(t)+z,t)
\right|_{z\downarrow -\xi}\right|\left| \left. \left.
\frac{d^{s-1-i}}{d\xi^{s-1-i}}
\right(\rho(-\xi_1)-\rho(-\xi_2)\right)\right|\\
&\leq C||\partial_t u||^{(2\beta)}|\xi_1-\xi_2|^{2\beta-[2\beta]},
\end{split}
\end{equation}
in which $C$ is a positive constant and $|| \cdot||^{(2\beta)}$ is
the H\"{o}lder norm (see Definition \ref{def:holder space}). On
the other hand, it follows from (\ref{eq:ind-assumption_2}) that
\begin{equation}\label{eq:est-int-ut-2}
\begin{split}
&\int_{-\xi_2}^{+\infty} \left| \left. \left.
\frac{\partial^s}{\partial \xi^s} \frac{\partial^r}{\partial
t^r}\frac{\partial}{\partial t}
\right(u(\xi_1+b(t)+z,t)-u(\xi_2+b(t)+z,t)\right) \right|\rho(z)dz
\\
&\leq ||\partial_t u||^{(2\beta)}
|\xi_1-\xi_2|^{2\beta-[2\beta]}\int_{-\xi_2}^{+\infty} \rho(z)dz
\leq ||\partial_t u||^{(2\beta)} |\xi_1-\xi_2|^{2\beta-[2\beta]}.
\end{split}
\end{equation}
Moreover, because $\rho(z)\in H^{\ell-1+2\alpha}(\R_-)$ for
$\ell\geq 1$ or $\int_{-\infty}^u \rho(z) dz \in
H^{2\alpha}(\R_-)$, we have $\int_{-\infty}^u \rho(z) dz \in
H^{\ell+2\alpha}(\R_-)$ for $\ell \geq 0$. In particular, using
$2\beta\leq 2\alpha+\ell$, we can see $\int_{-\infty}^u \rho(z) dz
\in H^{2\beta-[2\beta]}(\R_-)$. As a result,
\begin{equation}\label{eq:est-int-ut-3}
\begin{split}
&\int_{-\xi_1}^{-\xi_2} \left|\frac{\partial^s}{\partial \xi^s}
\frac{\partial^r}{\partial t^r}\frac{\partial}{\partial t}
u(\xi_1+b(t)+z,t)\right|\rho(z) dz \\
&\leq ||\partial_t u||^{(2\beta)} \left(\int_{-\infty}^{-\xi_2}
\rho(z)dz - \int_{-\infty}^{-\xi_1} \rho(z)dz \right) \leq
\tilde{C}||\partial_t u||^{(2\beta)}
|\xi_1-\xi_2|^{2\beta-[2\beta]},
\end{split}
\end{equation}
where $\tilde{C}$ is also a positive constant. Plugging the
estimates (\ref{eq:est-int-ut-1}) - (\ref{eq:est-int-ut-3}) into
(\ref{eq:est-int-ut}), we observe that (\ref{eq:holder-ut-xi})
holds.

When $\xi$ is fixed, using (\ref{eq:higher-der-ut}), it directly
follows from (\ref{eq:ind-assumption_1}) and (\ref{eq:ind-assumption_2}) that
\begin{equation}
 \frac{\partial^s}{\partial \xi^s}\frac{\partial^r}{\partial t^r}
 \int_{\R} \frac{\partial}{\partial t} u(\xi+b(t)+z,t) \nu(dz) \in
 H^{\beta-\frac{2r+s}{2}}([\delta',T]), \quad \text{for } 2\beta-2<2r+s<2\beta.
 \label{eq:holder-ut-t}
\end{equation}
Now, (\ref{eq:holder-ut-xi}) and (\ref{eq:holder-ut-t}) imply that
\begin{equation}
\int_{\R} \frac{\partial }{\partial t} u(\xi+b(t)+z,t)\nu(dz) \in
H^{2\beta,\beta}([0,+\infty)\times[\delta', T]).
\label{eq:holder-ut}
\end{equation}

Let $v(\xi,t)$ be a bounded solution of the boundary value problem
(\ref{eq:v_fix_1}) with the initial condition $v(\xi, \delta') = \partial_t u(\xi+b(\delta'), \delta')$. The uniqueness in Remark
\ref{remark:uniqueness} implies that
\begin{equation}\label{eq:v-eqs-toubt}
v(\xi,t)= \frac{\partial}{\partial t} u(\xi+b(t),t), \quad (\xi,
t) \in [0,+\infty) \times [\delta', T].
\end{equation}
As a result,  the assumption (\ref{eq:ind-assumption_2}) implies
that
\begin{equation}\label{eq:reg-v}
v(\xi,t)\in H^{2\beta,\beta}([0,+\infty)\times [\delta',T]).
\end{equation}
We will apply Lemma \ref{lemma:regular-improve} to
(\ref{eq:v_fix_1}) - (\ref{eq:v_fix_3}) with  $\phi(\xi+z,t) =
\partial_t u(\xi+b(t)+z,t) $, $a= \sigma^2/2$,
$\ell=-\left(\mu +b'(t)- \sigma^2/2\right)$, $c=r+\lambda$ and
$d=\lambda$. Thanks to (\ref{eq:ind-assumption_1}), the
coefficient $l$ belongs to $H^{\beta}([\delta,T])$. The other
coefficients already happen to reside there since they are
constants. Along with (\ref{eq:holder-ut}) and (\ref{eq:reg-v}),
Lemma \ref{lemma:regular-improve} yields
\begin{equation}\label{eq:reg-imp-v}
 v(\xi,t)\in
 H^{2+2\beta,1+\beta}([0,+\infty)\times[\delta'',T])\quad
 \text{for any } \delta''>\delta'>\delta,
\end{equation}
which implies that
\begin{equation}\label{eq:reg-imp-vxi}
\frac{\partial}{\partial \xi}v(0,t)\in
H^{\frac12+\beta}([\delta'',T]),
\end{equation}
and
\begin{equation}\label{eq:reg-imp-ut}
 \frac{\partial }{\partial t} u(\xi+b(t),t) \in
 H^{2+2\beta,1+\beta} ([0,+\infty)\times [\delta'',T]),
\end{equation}
by (\ref{eq:v-eqs-toubt}).

Using (\ref{eq:v_fix_4}) and (\ref{eq:reg-imp-vxi}), we will
improve the regularity of $b(t)$ in the following.  From
(\ref{eq:def-f-t-f}) we have
\begin{equation}\label{eq:f(b,t)}
\begin{split}
 f(b(t),t) &= \lambda
 \int_{\R}u(b(t)+z,t)\nu(dz) \\
 &= \lambda \int_0^{+\infty} u(b(t)+z,t)\nu(dz)
 + \lambda\int_{-\infty}^0
 (K-e^{b(t)+z})\nu(dz).
\end{split}
\end{equation}
Along with (\ref{eq:ind-assumption_1}) and
(\ref{eq:ind-assumption_3}), we can see from (\ref{eq:f(b,t)})
that
\begin{equation}\label{eq:reg-imp-f(b,t)}
f(b(t),t)\in H^{1+\beta}([\delta'',T]).
\end{equation}
Together with (\ref{eq:ind-assumption_1}), (\ref{eq:reg-imp-vxi})
and (\ref{eq:reg-imp-f(b,t)}), we can see from the identity
(\ref{eq:v_fix_4}) that $b'(t)\in H^{\frac12 +\beta}([\delta'',
T])$ for any $\delta''>\delta'$. It in turn implies that
\begin{equation}\label{eq:reg-imp-b}
 b(t)\in H^{\frac32+\beta}([\delta'',T]).
\end{equation}

{\bf Step 3.} Let us investigate $u(\xi+b(t),t)$. For any $r,s\geq
0$, $2r+s<2+2\beta$, we have
\begin{eqnarray*}
&&\frac{\partial^s}{\partial \xi^s} \frac{\partial^r}{\partial
t^r}
\int_{\R} u(\xi+b(t)+z,t) \nu(dz) \\
&&= \frac{\partial^s}{\partial \xi^s} \frac{\partial^r}{\partial
t^r} \int_{-\xi}^{+\infty} u(\xi+b(t)+z,t)\rho(z)dz +
\frac{\partial^s}{\partial \xi^s} \frac{\partial^r}{\partial t^r}
\int_{-\infty}^{-\xi}
u(\xi+b(t)+z,t)\rho(z)dz\\
&&= 1_{\{s\geq 1\}}\sum_{i=0}^{s-1} \left[
\frac{\partial^i}{\partial \xi^i}\left. \frac{\partial^r}{\partial
t^r} u(\xi+b(t)+z,t) \right|_{z\downarrow -\xi} -
\frac{\partial^i}{\partial \xi^i} \left.
\frac{\partial^r}{\partial t^r} u(\xi+b(t)+z,t) \right|_{z\uparrow
-\xi}\right] \frac{d^{s-1-i}}{d
\xi^{s-1-i}}\rho(-\xi)\\
&&\quad + \int_{-\xi}^{+\infty} \frac{\partial^s}{\partial \xi^s}
\frac{\partial^r}{\partial t^r} u(\xi+b(t)+z,t)\rho(z) dz +
\int_{-\infty}^{-\xi} \frac{\partial^s}{\partial \xi^s}
\frac{\partial^r}{\partial t^r} u(\xi+b(t)+z,t)\rho(z) dz,
\end{eqnarray*}
for any $\xi\geq 0$.
It is worth noticing that $\partial_{\xi}^i\partial_t^r
u(\xi+b(t)+z,t)|_{z\downarrow -\xi}\neq
\partial^i_{\xi}\partial^r_{t}u(\xi+b(t)+z,t)|_{z\uparrow -\xi}$ for some $i$ and
$r$. Following the same arguments that lead up to
(\ref{eq:holder-ut}), we can show
\begin{equation}
 \int_{\R} u(\xi+b(t)+z,t)\nu(dz) \in H^{2+2\beta, 1+\beta} ([0,+\infty)\times[\delta',
 T]), \label{eq:est-int-u}
\end{equation}
given $1+2\beta\leq 2\alpha+\ell -1$.

Now, we can apply Lemma \ref{lemma:regular-improve} to the
differential equation $u(\xi+b(t),t)$ satisfies, taking
(\ref{eq:ind-assumption_3}) and (\ref{eq:reg-imp-b}) into account.
This results in
\begin{equation}\label{eq:reg-imp-u}
u(\xi+b(t),t)\in
H^{3+2\beta,\frac32+\beta}([0,+\infty)\times[\delta''',T]),
\end{equation}
for any $\delta'''>\delta''$. As a result, we have improved the
regularities from (\ref{eq:ind-assumption_1}),
(\ref{eq:ind-assumption_2}) and (\ref{eq:ind-assumption_3}) to
(\ref{eq:reg-imp-b}), (\ref{eq:reg-imp-ut}) and
(\ref{eq:reg-imp-u}), respectively.

{\bf Step 4.} For any $\epsilon>0$, we apply Steps 2 and 3
inductively starting from $\beta=\alpha$ in Step 1.  Let $n$ be
the number of time we apply Steps 2 and 3. Let
$\delta_1'=\delta'$, in which $\delta'>0$ is as in Step 1. Running
Step 2 and 3 once, we obtain two constants $\delta_1''$ and
$\delta'''_1$ such that (\ref{eq:reg-imp-b}), (\ref{eq:reg-imp-u})
hold with $\beta=\alpha$.  In the $n$-th time, $n\geq 2$, we
choose $\delta_n' = \delta_{n-1}'''$ and  $\delta_n''' >
\delta_n'' >\delta_n'$, such that $\delta_n'''<\epsilon$ for any
$n$ so that $[\epsilon, T]\subset [\delta'''_n,T]$.

The application of Step 2 for the $n$-th time will give us that
$b(t)\in H^{1+\alpha+\frac{n}{2}}([\epsilon, T])$. Applying Step 2
for $\ell+1$ and Step 3 for $\ell$ times the result follows.
\end{proof}

\begin{remark}\label{remark:higher-reg-u}
\begin{itemize}
\item[(i)] The previous proof has also shown the higher order
regularity of $u(x,t)$, i.e. $u(\xi+b(t),t)\in H^{2+2\alpha+\ell,
1+\alpha+\frac{\ell}{2}}([0,+\infty)\times [\epsilon, T])$, for
any $\epsilon>0$, under the assumptions of Theorem
\ref{theorem:higher-reg-b}, . \item[(ii)] Note that $b(t)\in
C^1((0,T])$ without any assumption on the density $\rho(z)$. If
$\rho(z)\in H^{2m-1+2\alpha}(\R_-)$ for some $m\geq 1$, then
$b(t)\in H^{\frac32 + m+ \alpha}([\epsilon, T])$. From Definition
\ref{def:holder space} and the arbitrary choice of $\epsilon$, we
have that $b(t)\in C^{m+1}((0,T])$ under this assumption.
\end{itemize}

\end{remark}

As a corollary of Theorem \ref{theorem:higher-reg-b}, we have the
following sufficient condition for the infinitely
differentiability of $b(t)$.

\begin{corollary}\label{cor:inf-dif-b}
Let $b(t)$ be the free boundary in (\ref{eq:u_1}) -
(\ref{eq:u_4}).  Assume that $\nu$ has a density, i.e. $\nu(dz)=
\rho(z) dz$. If $\rho(z)\in C^{\infty}(\R_-)$ with
$\frac{d^{\ell}}{dz^{\ell}} \rho(z)$ bounded for each $\ell \geq
1$, but not necessarily uniformly, then $b(t)\in
C^{\infty}((0,T])$.
\end{corollary}

\begin{proof}
For any $m \geq 1$ with $\rho(z)\in C^{2m+1}(\R_-)$ and
derivatives of $\rho(z)$ up to order $2m+1$ are bounded, it
follows from Definition \ref{def:holder space} that $\rho(z)\in
H^{2m-1+2\alpha}(\R_-)$. As a result of Remark
\ref{remark:higher-reg-u} (ii), we have $b(t)\in C^{m+1}((0,T])$.
\end{proof}

\begin{remark}\label{remark:Merton-Kou}
There are two well-known examples of jump diffusion models in the
literature, Kou's model and Merton's model (see \cite{cont},
p.111), in which the density $\rho(z)$ is double exponential and
normal, respectively. For both of these densities, it is easy to
see that the conditions for Corollary \ref{cor:inf-dif-b} are
satisfied. Therefore, the free boundaries in both models are
infinitely differentiable.
\end{remark}

\section{The boundaries of the approximating free boundary problems introduced by Bayraktar [2008]}
In this section, we want to show that the
approximating free boundaries $b_n(t)$, constructed in
\cite{bayraktar-finite-horizon}, have regularity
properties similar to the free boundary $b(t)$.

\cite{bayraktar-finite-horizon} constructed a monotone increasing
sequence $\{u_n\}_{n\geq 0}$ that converges to the unique solution
$u(x,t)$ of the parabolic integro-differential equation
(\ref{eq:u_1}) - (\ref{eq:u_4}), uniformly. In this sequence,
$u_0(x,t)=(K-e^x)^+$, and each $u_{n}(x,t)$ $(n\geq 1)$ is the
unique classical solution of the following parabolic differential
equation:
\begin{eqnarray}
 &&\L_{\D}u_n \triangleq \frac{\partial u_{n}}{\partial t} -\frac12 \sigma^2 \frac{\partial^2 u_{n}}{\partial
 x^2} -\left(\mu-\frac12 \sigma^2\right) \frac{\partial
 u_{n}}{\partial x} + (r+\lambda)u_{n} =
 f_n(x,t), \quad x>b_{n}(t), \label{eq:u_n_1}\\
 && u_{n}(b_{n}(t),t) = K- e^{b_{n}(t)}, \quad
 t\in(0,T], \label{eq:u_n_2}\\
 && u_n(x,0) =(K-e^x)^+, \quad x\geq b_n(0), \label{eq:u_n_4}
\end{eqnarray}
in which
\begin{equation}\label{eq:def-fn}
f_n(x,t) \triangleq \lambda \int_{\R} u_{n-1}(x+z,t) \nu(dz),
\end{equation}
and the free boundary $b_{n}(t) \triangleq
\log\left(s_{n}\left(T-t\right)\right)$ is defined in terms of
$s_{n}(\cdot)$, which is the approximating free boundary in
\cite{bayraktar-finite-horizon}. Moreover, the smooth fit property
is also satisfied for each $u_n$, i.e.
\begin{equation}\label{eq:u_n_3}
 \frac{\partial}{\partial x} u_{n}(b_{n}(t),t) =
 -e^{b_{n}(t)},\quad t\in(0,T].
\end{equation}
In the region $\{(x,t) |\, x<b_n(t), t\in(0,T]\}$, one also has
that
\begin{equation}\label{eq:u_n_5}
 \L_{\D} u_n (x,t) - f_n (x,t) \geq 0.
\end{equation}
We can define the approximating continuation regions
$\mathcal{C}_n$ and the stopping regions $\mathcal{D}_n$ as
follows
\begin{eqnarray*}
\mathcal{C}_n \triangleq \{(x,t)\,|\, b_n(t)<x<+\infty, 0<t\leq
T\}, \quad \mathcal{D}_n\triangleq  \{(x,t)| -\infty<x\leq b_n(t),
0<t\leq T\}, \quad \text{ for all } n\geq 1.
\end{eqnarray*}
Since $\{u_n\}_{n\geq 0}$ is a monotone increasing sequence, the
approximating free boundary $\{b_n\}_{n\geq 1}$ is a monotone
decreasing sequence. As a result, we have $\cup_{n\geq 1}
\mathcal{C}_n = \mathcal{C}$ and $\cap_{n\geq 1} \mathcal{D}_n =
\mathcal{D}$.

The approximating sequences $\{u_n\}_{n\geq 1}$ and
$\{b_n\}_{n\geq 1}$ have the similar properties with the value
function $u$ and its free boundary $b$. Proposition~\ref{prop:u_t
positive}, Lemmas~\ref{lemma:b_strictly}, \ref{lemma:lip-holder}
and \ref{lemma:u_t_inf} have their analogous versions for $u_n$
and $b_n$ via the same proofs only replacing the integral term $f$
in \eqref{eq:u_1_para} by $f_n$ in \eqref{eq:def-fn}. 
Proposition~\ref{prop:u_t} and Lemma~\ref{lemma:ut-bounded}, on the other hand, can be
slightly modified as follows:

\begin{proposition}\label{prop:un_t} For all $n\geq 1$,\\
{\bf (i)} If $\partial_t u_{n-1}(x,t)$ is bounded in
$\R\times[\epsilon, T]$ for any $\epsilon>0$, then $\partial_t
u_n(x,t)$ is continuous in $\R\times(0,T]$ and
\begin{eqnarray}
 \lim_{x\downarrow b_{n}(t)} \frac{\partial}{\partial t}
 u_{n}(x,t)&=&0  \label{eq:u_n_t}.
\end{eqnarray}
{\bf (ii)} On the other hand, if $\lim_{x\downarrow b_n(t)}\partial_t u_n(x,t) =0$
for $t\in(0,T]$ and $\partial_t
u_n(x,t)$ is continuous in $\R\times(0,T]$, then $\partial_t u_n(x,t)$ is uniformly bounded in
$\R\times[\epsilon, T]$, for any $\epsilon>0$.
\end{proposition}
\begin{proof}
See Appendix~\ref{app:A.3} for the proof of (i). Under the
assumption that $\lim_{x\downarrow b_n(t)} \partial_tu_n(x,t) =0$ for
$t\in(0,T]$, we have $\partial_t u(x,t)$ is bounded in the domain
$\{(x,t)\, | \, b_n(t)\leq x\leq X_0, \epsilon \leq t\leq T\}$ for
any $\epsilon\geq 0$ and $X_0>\log K$. Then the rest of the proof of (ii) is similar to
the proof of Lemma~\ref{lemma:ut-bounded}.
\end{proof}
\begin{remark}
To show that assumptions in both (i) and (ii) are satisfied for
all $u_n$, $n\geq 1$, we need to walk through (i) and (ii)
successively. Starting from $\partial_t u_0(x,t) =0$ (since
$u_0(x,t)=(K-e^x)^+$), (i) tells us that $\lim_{x\downarrow
b_1(t)}\partial_t u_1(x,t) =0$ and $\partial_t
u_1(x,t)$ is continuous in $\R\times(0,T]$. Then it follows from (ii) that
$\partial_t u_1(x,t)$ is bounded in $\R\times[\epsilon, T]$ for
any $\epsilon>0$. This result feeds back to (i). Now, as a result of an induction argument
it can be seen that assumptions in both (i) and (ii) are satisfied for all $n$.
\end{remark}

Results similar to Lemmas \ref{lemma:j},
 \ref{lemma:level-B} and Corollary \ref{cor:u_xx_pos} can also be
 shown to hold for each $u_n$, $n\geq 1$. Defining
 \begin{equation*}
 \begin{split}
  J_n(x,t) &\triangleq qe^x - rK + \lambda \int_{\R} \left[u_{n-1}(x+z,t) + e^{x+z} -
  K\right]\nu(dz), \quad x\in \R, t\in[0,T],\\
  B_n(t) & \triangleq \left\{x : J_n(x,t)=0, t\in [0,T]\right\}.
 \end{split}
 \end{equation*}
 we obtain the following:
 \begin{eqnarray}
  && \L_{\D}u_n(x,t)- \lambda \int_{\R}
  u_{n-1}(x+z,t) \nu(dz) = -J_n(x,t), \quad
  x<b_n(t), t\in[0,T],\label{eq:lun=jn}\\
  && x\rightarrow J_n(x,t) \text{ is strictly increasing and }  t\rightarrow J_n(x,t) \text{ is non-decreasing for } (x,t)\in
  \R\times[0,T], \label{eq:jn-inc}\\
  && B_n(t) >b_n(t), \quad t\in(0,T], \\
  && \lim_{x\downarrow b_n(t)} \frac{\partial^2}{\partial x^2}
  u_n(x,t) > - e^{b_n(t)}, \quad t\in(0,T]. \label{eq:un_xx_pos}
 \end{eqnarray}

Moreover, as we can see in the following Proposition, the
approximating free boundaries $b_n$ have the same critical value
as $b$ at 0.

\begin{proposition}\label{prop:bn(0)}
For the approximating sequence $b_n(t)$, we have
\begin{equation}
b_n(0+) \triangleq \lim_{t\rightarrow 0^+}b_n(t) =
\min\{\log{K},B(0)\} = \left\{\begin{array}{ll} \log{K}, & r\geq
q+\lambda\int_{\R_+}(e^z-1)\nu(dz) \\
B(0), & r< q+\lambda\int_{\R_+}(e^z-1)\nu(dz)
\end{array}\right. , \label{eq:b-n(0)}
\end{equation}
in which $B(0)$ the unique solution of \eqref{eq:B(0)-eq}.
\end{proposition}
\begin{proof}
When $x<b_n(t) (t> 0)$, it follows from \eqref{eq:u_n_5},
\eqref{eq:lun=jn} and \eqref{eq:jn-inc} that
\begin{equation*}
 0 \leq \L_{D} u_n(x,t) -\lambda \int_{\R}
 u_{n-1}(x+z,t) \nu(dz) = - J_n(x,t) \leq -
 J_n(x,0) = - J_0(x).
\end{equation*}
The fact that $J_0(B(0))=0$ and $x\rightarrow J_0(x)$ is strictly
increasing tells us that $x \leq B(0)$. Hence $b_n(t)\leq B(0)$ thanks to the
choice of $x$. It is also clear that $b_n(t)\leq \log{K}$. Then we
obtain
\begin{equation}\label{eq:b_n(0)}
 b_n(0+) \leq \min\{\log{K}, B(0)\}.
\end{equation}
 Now, the
corollary results from combining (\ref{eq:b(0)}) and
(\ref{eq:b_n(0)}), since $\{b_n\}_{n\geq 1}$ is a
decreasing sequence of functions.
\end{proof}

Furthermore, the H\"{o}lder continuity in
Theorem~\ref{theorem:holder-cont} also holds for $b_n$, $n\geq 1$.
In the proof of Lemma~\ref{lemma:b_holder_est}, we only need to
replace $c$ in \eqref{eq:chs-c} by $\min\left\{-2/{\sigma^2}
J_n(x,t) | \, b_n(t)<x<B_n(t), \epsilon \leq t \leq T \right\}>0$.
On the other hand, results in Lemma~\ref{lemma:u_xt cont} also
hold for $\partial_{xt} u_n$, $n\geq 1$. Therefore, combining with
\eqref{eq:un_xx_pos}, we have from \eqref{eq:u_n_3} that
\begin{proposition}
 $b_n(t)\in C^1(0,T]$, $n\geq 1$.
\end{proposition}

Finally, using the following representation
\begin{equation}\label{eq:def-bn'}
 b_n'(t)= - \frac{\frac{\sigma^2}{2}\frac{\partial^2}{\partial
 x \partial t}u_n(b_n(t)+,t)}{(\mu-r-\lambda)e^{b_n(t)}+ (r+\lambda)K -
 f_n(b_n(t),t)}, \quad t\in(0,T],
\end{equation}
one can follow the proof of Lemma~\ref{lemma:b'-holder} to show
that there is $\alpha \in (0,1/2)$ such that
\[
 b_n(t)\in H^{1+\alpha}([\delta,T]), \quad \text{ for any } \delta>0.
\]

\appendix
\renewcommand{\theequation}{A-\arabic{equation}}
\renewcommand{\thetheorem}{A-\arabic{theorem}}
\renewcommand{\thedefinition}{A-\arabic{definition}}
\renewcommand{\thelemma}{A-\arabic{lemma}}
\renewcommand{\theremark}{A-\arabic{remark}}
\setcounter{section}{0}%
\renewcommand{\thesection}{\Alph{section}}%
\section{}

\subsection{Proof of Lemmas {\ref{lemma:lip-holder}}, {\ref{lemma:ut-bounded}} and {\ref{lemma:u_t_inf}}}\label{app:A.1}
\begin{proof}[{\bf Proof of Lemma {\ref{lemma:lip-holder}}}]
 The inequality (\ref{eq:holder-u}) is clear, because we have
 \[
  |u(x,t)-u(x,s)|= \left|V(e^x, T-t)- V(e^x,
  T-s)\right|\leq
  D |t-s|^{\frac12}.
 \]

 In order to prove (\ref{eq:lip-u}), it suffices to check that $\partial_x
 u(x,t)$ is uniformly bounded in the domain $\R\times[0,T]$. Choose a constant
 $X>\log{K}+1$, we will first prove $\partial_x u(x,t)$ is uniformly
 bounded in $[X,+\infty)\times[0,T]$. Let us consider a cut-off
 function $\eta(x)\in C^{\infty}(\R)$, such that $\eta(x)=0$ when
 $x\leq X-1$ and $\eta(x)=1$ when $x\geq X$. Using
 (\ref{eq:u_1_para}) we see that $v(x,t)=\eta(x)u(x,t)$ satisfies
 \begin{eqnarray*}
  &&\L_{\D} v = \eta(x)f(x,t) + \tilde{f}(x,t), \\
  && v(x,0)= \eta(x)(K-e^x)^+,
 \end{eqnarray*}
 where
 \begin{equation}\label{eq:def-f-t-f}
  f(x,t)= \lambda \int_{\R} u(x+z,t)\nu(dz),
  \quad \tilde{f}(x,t)= - \frac12 \sigma^2 \left(\eta'' u+ 2\eta' \frac{\partial u}{\partial x}\right)
  - \left(\mu - \frac12 \sigma^2\right)\eta' u.
 \end{equation}
 It is worth noticing that the term $\eta' \partial_x u$ in the expression for $\tilde{f}$
 vanishes outside a compact domain. Since we also have that  $u(x,t)\leq
 K$, both $f(x,t)$ and
 $\tilde{f}(x,t)$ are bounded in $\R\times[0,T]$.

 Let $G(x,t;y,s)$ be the Green function corresponding to the differential
 operator $\L_{\D}$. We can represent $v(x,t)$ in terms of $G$ as
 \begin{equation}\label{eq:rep-v}
  v(x,t)= \int_{\R} dy\, G(x,t;y,0) \eta(y)(K-e^y)^+ + \int_0^t
  ds\int_{\R} dy\, G(x,t;y,s)\left(f(y,s) \eta(y)+ \tilde{f}(y,s)\right).
 \end{equation}
 The first term on the right-hand-side of (\ref{eq:rep-v}) will vanish by
 the choice of $\eta(y)$. On the other hand, Green function
 $G(x,t;y,s)$ satisfies
 \[
  \left|\partial_x G(x,t;y,s)\right| \leq c(t-s)^{-1}
  \exp\left(-c \frac{|x-y|^2}{t-s}\right),
 \]
 for some positive constant $c$, (see Theorem 16.3 in page 413 of
 \cite{lad}). Since $\int_{\R} dy\, \exp(-c\frac{(x-y)^2}{t-s}) \leq d\, (t-s)^{\frac12}$ for some other positive constant $d$, we have that
 \[
  \int_0^t ds \int_{\R} dy\, |\partial_x G(x,t;y,s)| \leq \int_0^t
  ds\, \tilde{c}(t-s)^{-\frac12} = 2\tilde{c}\,t^{\frac12},
 \]
 Using this
 estimate and the boundness of $f$ and $\tilde{f}$, the Dominated
 Convergence Theorem implies that
 \[
  \partial_x v(x,t) = \int_0^t ds\int_{\R} dy\, \partial_x
  G(x,t;y,s)(f(y,s) \eta(y)+\tilde{f}(y,s)),
 \]
 which is uniformly bounded. On the other hand, $\partial_x v= \eta' u + \eta \partial_x
 u$. By our choice of $\eta(x)$, we  have that $\partial_x u(x,t)$ is
 uniformly bounded on $[X,+\infty)\times[0,T]$.

 Moreover, in the stopping region $\mathcal{D}$, we have $\partial_x u(x,t)= -
 e^x$. This implies that $0>\partial_x u(x,t)\geq -e^{b(t)} \geq -K$. On the other hand,  since it is continuous $\partial_x u$ is also bounded in
 the compact closed domain $\{(x,t)|b(t)\leq x\leq X, 0\leq t\leq
 T\}$. As a result
 we have that $\partial_x u(x,t)$ is uniformly bounded in
 $\R\times[0,T]$.
\end{proof}

\begin{proof}[{\bf Proof of Lemma \ref{lemma:ut-bounded}}]
Let us choose $X_0$ such that $X_0 > \log{K}$. We will first prove
that $\partial_t u(x,t)$ is uniformly bounded in the domain
$[X_0,+\infty)\times [0,T]$. Let $k(x,t)\in C^{\infty}_0(\R\times
[0,T])$ be such that
\[
 \partial_x k(x,t) |_{x=X_0} = \partial_x u(x,t) |_{x=X_0}, \quad t\in [0,T],
\]
and that $k(x,0)=0$, $x \in \mathbb{R}$. These two conditions on
$k$ are consistent since
 $\partial_x u(x,0)|_{x=X_0} =0$. The function
$v(x,t) \triangleq u(x,t)-k(x,t)$ satisfies
\begin{equation}\label{eq:boundary-v}
\partial_x v(x,t)|_{x=X_0}=0,
\end{equation}
and
\begin{equation}
 \L_{\D} v(x,t) = f(x,t) + g(x,t), \quad x>b(t), t\in (0,T],
\end{equation}
in which $g(x,t)=-\L_{\D}k(x,t)$ and $f$ is given by
(\ref{eq:def-f-t-f}). Let us define the even extension of $v(x,t)$
with respect to the line $x=X_0$ as
\begin{equation}\label{eq:v-evenext}
 \hat{v}(x,t) \triangleq \left\{\begin{array}{ll} v(x,t) & x\geq X_0, \\ v(2X_0-x,t) & x<X_0.
 \end{array}\right.
\end{equation}
We similarly define $\hat{f}(x,t)$ and $\hat{g}(x)$. From
(\ref{eq:boundary-v}) and (\ref{eq:v-evenext}), we have
$\hat{v}(x,t)\in C^{2,1}(\R\times(0,T])$ and that it satisfies the
equation
\begin{eqnarray*}
&&\L_{\D}\hat{v} = \hat{f}(x,t)+ \hat{g}(x,t), \quad (x,t)\in
\R\times
(0,T],\\
&& \hat{v}(x,0) = 0, \quad x\in \R.
\end{eqnarray*}
Here the initial condition follows from (\ref{eq:u_4}) and the
choice of $X_0$ and $k(x,t)$.

It follows from (\ref{eq:lip-u}) and (\ref{eq:holder-u}) that
$f(x,t)$ is uniformly Lipschitz in $x$ and semi-H\"{o}lder
coninuous in $t$. So for any $x_1<x_2$, if we have either $x_2\leq
X_0$ or $X_0\leq x_1$, then
\[
 \left|\hat{f}(x_1,t)-\hat{f}(x_2,t)\right| \leq \lambda C (x_2-x_1),
\]
for the same constant $C$ as in (\ref{eq:lip-u}). On the other
hand, if $x_1<X_0<x_2$, then
\begin{eqnarray*}
|\hat{f}(x_1,t)- \hat{f}(x_2,t)| &\leq&
|\hat{f}(x_1,t)-\hat{f}(X_0,t)| +
|\hat{f}(X_0,t)-\hat{f}(x_2,t)| \\
&\leq& \lambda C(X_0-x_1) + \lambda C(x_2-X_0) = \lambda
C(x_2-x_1).
\end{eqnarray*}
As a result of the last two equations we observe that $\hat{f}(x,
t)$ is uniformly Lipschitz in its first variable. It is also clear
that $\hat{f}(x, t)$ is semi-H\"{o}lder continuous in its second
variable. Thus, it follows from Definition \ref{def:holder space}
that
\[
 \hat{f}(x,t)\in H^{\alpha,\frac{\alpha}{2}}(\R\times[0,T]),
 \quad \text{ for some } 0<\alpha<1.
\]
On the other hand, $\hat{g}(x,t)\in
H^{\alpha,\alpha/2}(\R\times[0,T])$, because $k(x,t)\in
C^{\infty}_0 (\R\times[0,T])$. Combining with the assumption
\eqref{eq:ass-sigma} on $\sigma$, the regularity property of
parabolic differential equation (see Theorem 5.1 in page 320 of
\cite{lad}) implies that
\[
 \hat{v}(x,t) \in H^{2+\alpha,
 1+\frac{\alpha}{2}}(\R\times[0,T]).
\]
In particular, $u(x,t)\in H^{2+\alpha,
1+\alpha/2}([X_0,+\infty)\times[0,T])$. As a result,  in
$[X_0,+\infty)\times[0,T]$, $\partial_t u(x,t)$ is uniformly
bounded by the H\"{o}lder norm of $u(x,t)$ . Now, the result
follows from the continuity of $\partial_t u(x,t)$ inside domain
$\{(x,t)\,| \,b(t)\leq x\leq X_0, \epsilon\leq t\leq T\}$ for any
$\epsilon >0$ (see Proposition \ref{prop:u_t}).
\end{proof}

\begin{proof}[{\bf Proof of Lemma \ref{lemma:u_t_inf}}]
Let $X_0 > \log{K}$  be the same as in the proof of Lemma
\ref{lemma:ut-bounded}, again choose a cut-off function
$\eta(x)\in C^{\infty}(\R)$, such that $\eta(x)=1$ when $x\geq
2X_0$ and $\eta(x)=0$ when $x\leq X_0$. Then formally the function
$\eta(x)\partial_t u(x,t)$ satisfies the following Cauchy problem
\[
 \L_{\D} w = \eta(x)h(x,t) + \tilde{h}(x,t), \quad (x,t)\in \R\times [t_0,T],
\]
where
\[
 h(x,t)= \lambda\int_{\R}\partial_tu(x+z,t)\nu(dz), \quad \tilde{h}(x,t)= -\frac12 \sigma^2 \left(2\eta'
 \partial_x\partial_t
u+ \eta''\partial_t u\right)- \left(\mu-\frac12
\sigma^2\right)\eta'\partial_tu,
\]
and we choose $\eta(x)\partial_t u(x,t_0)$, for some $t_0 \in
(0,T)$, as the initial condition. It follows from Theorem 3.1 in
page 346 of \cite{garr-mena-1} that this Cauchy problem has an
unique classical solution, we call it $w$. On the other hand, we
have $w(x,t) = \eta(x) \partial_t u(x,t)$. Indeed, it is easy to
check that $\int_{t_0}^t w(x,s)ds$ is the unique classical
solution of the Cauchy problem
\[
 \L_{D} v = \int_{t_0}^t ds\, \left(\eta(x)h(x,s) +
 \tilde{h}(x,s)\right) + \eta(x)\partial_t u(x,t_0), \quad
 v(x,t_0)=0.
\]
Note that $\eta(x)\left[u(x,t)-u(x,t_0)\right]$ is another
classical solution. Therefore $w(x,t)= \eta(x)\partial_t u(x,t)$
by the uniqueness.

Using the Green function $G(x,t;y,s)$ corresponding to the
differential operator $\L_{\D}$, the solution $w(x,t)$ can be
represented as
\begin{equation}\label{eq:rep-for-w}
 w(x,t) = \int_{\R}dy\, G(x,t;y,t_0)w(y,t_0) + \int_{t_0}^{t}
 ds\int_{\R} dy\, G(x,t;y,s) (\eta(y)h(y,s)+\tilde{h}(y,s)),
\end{equation}
for all $(x,t)\in \R\times(t_0,T]$. Since the Green function
satisfies
\[
 |G(x,t;y,s)|\leq
 C(t-s)^{-\frac12}\exp\left(-\frac{c(x-y)^2}{t-s}\right), \quad
 (y,s)\in \R\times[0,t).
\]
The first term in (\ref{eq:rep-for-w}) is bounded, as long as
$w(y,t_0)$ is uniformly bounded. The contribution of
$\eta'\partial_x\partial_t u$ (in the expression for $\tilde{h}$)
to $w$ is given by,
\[
-\int_{\R}dy \,G(x,t;y,s) \eta'(y) \frac{\partial^2}{\partial
y\partial s} u(y,s)=\int_{\R} dy \frac{\partial }{\partial
y}[G(x,t;y,s) \eta'(y)]\frac{\partial }{\partial s} u(y,s).
\]
Now it follows from Lemma \ref{lemma:ut-bounded} that both
$w(x,t_0)$ and $h(x,t)$ are uniformly bounded for $x\in \R, t\in
[t_0,T]$. We also have that $\eta'$ and $\eta''$ vanish outside
$[X_0,2X_0]$. Since $\lim_{x\rightarrow +\infty} G(x,t;y,s)=0$ and
it can easily be shown that $ \lim_{x\rightarrow +\infty}
\partial_y G(s,t;y,s)=0$, the Dominated Convergence Theorem
implies that
\[
 \lim_{x\rightarrow +\infty} w(x,t) =0, \quad t\in (t_0,T].
\]
Then the statement follows from the choice of $\eta$.
\end{proof}

\subsection{Proof of Lemma {\ref{lemma:u_xt cont}}}\label{app:A.2}
We will first establish a one to one correspondence between
solutions of (\ref{eq:w_1}) and solutions of an integral equation
of Volterra type.

\begin{lemma}\label{lemma:correspondence}
(i) Let $G(x,t;y,s)$ be the Green function associated to the differential operator $\L_{\D}$ and
let us consider the following nonlinear integral equation of Volterra type,
\begin{equation}\label{eq:volterra}
 \left(1+\frac14 \sigma^2(b(t),t)\right) v(t) = -\int_{t_0}^t ds \,v(s) \frac12 \sigma^2(b(s), s)\, \partial_x G(b(t), t; b(s), s) +
 \sum_{i=1}^2 N_i(t), \quad t_0\leq t\leq T,
\end{equation}
where $N_1(t)= \int_{b(t)}^{+\infty} dy\, \partial_x
G(b(t),t;y,t_0) w(y,t_0)$ and $N_2(t)=\int_{t_0}^tds
\int_{b(s)}^{+\infty} dy\,
 \partial_x G(b(t),t;y,s)h(y,s)$.
 There exists a unique solution $v$ to (\ref{eq:volterra}). The function $v(t)$ is continuous.

(ii)  Let $w(x,t)$ be a classical solution of 
(\ref{eq:w_1}) on $[t_0, T]$ with the initial condition $w(x,t_0)=
\partial_t u(x,t_0)$, such that $t \rightarrow \partial_x w(b(t)+,
t)$ is continuous. Then there is a one to one correspondence
between $w(x,t)$ and $v(t)$. Moreover $\partial_x w(b(t)+,t)=
v(t)$, $t_0\leq t\leq T$.
\end{lemma}
The initial value of (\ref{eq:w_1}) may not be smooth.
This is the reason we take $w(x,t_0)=\partial_t u(x,t_0)$,
$0<t_0<T$, as the initial condition of (\ref{eq:w_1}) and consider
the differential equation on $t\in[t_0,T]$.

\begin{remark}
The correspondence in Lemma~\ref{lemma:correspondence} is well known for the Stefan problem on heat equation with Lipschitz continuous free boundary (see Section 1
Chapter 8 of \cite{friedman-1}). Along Friedman's line of proof,
we will extend the correspondence to our parabolic differential
equation with H\"{o}lder continuous free boundary.
\end{remark}

\noindent \emph{Proof of Lemma~\ref{lemma:correspondence}}.
\noindent Proof of (i). First, because $G(b(t), t; b(s), s)$ and
$\sigma(b(s),s)$ are continuous for $s \in(0,t)$ (see
\eqref{eq:ass-sigma}), it follows from the classical result on
Volterra equations (see \cite{rust}) that the integral equation
(\ref{eq:volterra}) has a unique solution $v(t)$ and it is
continuous with respect to $t\in[t_0,T]$, as long as $N_i(t)$,
$i=1,2$, are continuous with respect to $t$. It is not hard to
show these functions are indeed continuous, using the continuity
of $b(t)$ and the following estimates on the Green function $G$
and its derivatives: {\small
\begin{eqnarray*}
 &&|\partial_x^{\ell}G(x,t;y,s)| \leq C(t-s)^{-\frac{1+\ell}{2}}
 \exp\left(-c\frac{|x-y|^2}{t-s}\right),\\
 &&|\partial_x G(x,t;y,s)-\partial_x G(x,\tilde{t};y,s)| \leq
 C(t-\tilde{t})^{\frac{\alpha}{2}}(\tilde{t}-s)^{-\frac{2+\alpha}{2}}\exp\left(-c\frac{|x-y|^2}{t-s}\right),\\
 && |\partial_x G(x,t;y,s) - \partial_{\tilde{x}}
 G(\tilde{x},t;y,s)| \leq
 C|x-\tilde{x}|^{\alpha}(t-s)^{-\frac{2+\alpha}{2}}\exp\left(-c\frac{|x''-y|^2}{t-s}\right),
\end{eqnarray*}}
where $\ell =0,1$, $s<\tilde{t}<t$,
$|x''-y|=|x-y|\wedge|\tilde{x}-y|$, $0<\alpha<1$, $C$ and $c$ are
positive constants. These estimates are from Theorem 16.3 in page 413 of \cite{lad}. \\

\noindent Proof of (ii) Let us assume that  $w(x,t)$ is a
classical solution of (\ref{eq:w_1}).  As a result,  the following
Green's identity (see page 27 of \cite{friedman-1}) is satisfied
 \begin{equation}\label{eq:Green-identity}
 \begin{split}
 \frac{\partial}{\partial y}\left(\frac12 \sigma^2(y,s) G(x,t;y,s) \frac{\partial}{\partial y}w(y,s) - \frac12 \sigma^2(y,s) w(y,s)\frac{\partial}{\partial
 y}G(x,t;y,s) - w(y,s)G(x,t;y,s) \sigma \sigma_y(y,s)\right) \\
 - \frac{\partial}{\partial s}\left(G(x,t;y,s)w(y,s)\right)
 + \frac{\partial}{\partial y}\left(\left(\mu-\frac12 \sigma^2(y,s)\right)G(x,t;y,s)w(y,s)\right) =
 -G(x,t;y,s)h(y,s),
\end{split}
 \end{equation}
 where $\quad t_0 \leq s < t \leq T$, $x> b(t)$ and  $y>b(s)$.
Integrating both hand side of (\ref{eq:Green-identity}) over the
domain $b(s)< y< +\infty$, $t_0 < s< t-\epsilon$, we obtain
\begin{equation}\label{eq:int-Green}
\begin{split}
 &\int_{t_0}^{t-\epsilon} ds \lim_{y\rightarrow +\infty} \frac12
 \sigma^2(y,s)\, \partial_y w(y,s) \, G(x,t;y,s) - \int_{t_0}^{t-\epsilon}
 ds \,\frac12 \sigma^2(b(s), s)\,
 \partial_y w(b(s)+,s)\, G(x,t;b(s), s) \\
 &-\int_{t_0}^{t-\epsilon} ds \lim_{y\rightarrow +\infty} \frac12
 \sigma^2(y,s)\,
 w(y,s)\, \partial_y G(x,t;y,s) + \int_{t_0}^{t-\epsilon}
 ds \,\frac12 \sigma^2(b(s),s)\, w(b(s),s)\, \partial_y G(x,t;b(s),s)\\
 &-\int_{t_0}^{t-\epsilon} ds \lim_{y\rightarrow +\infty} w(y,s)
 \, G(x,t; y,s) \, \sigma \sigma_y(y,s) +
 \int_{t_0}^{t-\epsilon}ds \, w(b(s), s) \, G(x, t; b(s), s)\,
 \sigma \sigma_y(b(s),s)\\
 &- \int_{b(t-\epsilon)}^{+\infty} dy \left[G(x,t;y,t-\epsilon)w(y,t-\epsilon) -
 G(x,t;y,t_0)w(y,t_0)\right]\\
 &+ \int_{t_0}^{t-\epsilon}ds \left[\lim_{y\rightarrow +\infty}\left(\mu-\frac12 \sigma^2(y,s)\right)w(y,s)\,G(x,t;y,s)
 - \left(\mu-\frac12 \sigma^2(b(s),s)\right)
 w(b(s),s)\,G(x,t;b(s),s)\right]\\
 &= -\int_{t_0}^{t-\epsilon}ds \int_{b(s)}^{+\infty} dy\,
 G(x,t;y,s)h(y,s).
 \end{split}
 \end{equation}
 In the seventh term on the left of (\ref{eq:int-Green}), we used $w(x,t)=0$ when $x< b(t)$.
 Using the boundary and initial conditions for $w(x,t)$ and the facts that $\lim_{y\rightarrow +\infty} G(x,t;y,s)=0$
 and $\lim_{y\rightarrow +\infty} \partial_y G(x,t;y,s)=0$, letting $\epsilon \rightarrow 0$, we can write
\begin{equation}\label{eq:u(x,t)}
\begin{split}
 w(x,t) &= -\int_{t_0}^t ds\, \partial_x
 w(b(s)+,s)\, \frac12 \sigma^2(b(s),s)\,
 G(x,t;b(s),s) + \int_{b(t)}^{+\infty} dy\, G(x,t;y,t_0)
 w(y,t_0)\\
 & \hspace{0.4cm} +\int_{t_0}^{t}ds \int_{b(s)}^{+\infty} dy\,
 G(x,t;y,s)h(y,s) \\
 &\triangleq -M_0(x,t) + M_1(x,t) + M_2(x,t).
\end{split}
\end{equation}

Before differentiating both sides of (\ref{eq:u(x,t)}) with
respect to $x$, let us recall the jump identity: if $\rho(t)$,
$t_0\leq t\leq T$, is a continuous function and $b(t)$ is the
H\"{older} continuous with H\"{o}lder exponent $\alpha>\frac12$,
then for every $t_0\leq t \leq T$,
\begin{equation}\label{eq:jump identity}
\lim_{x\downarrow b(t)} \frac{\partial}{\partial x} \int_{t_0}^t ds\, \rho(s) G(x,t; b(s), s) = \frac12 \rho(t) + \int_{t_0}^t ds\, \rho(s) \left. \partial_x G(x,t; b(s), s) \right|_{x=b(t)}.
\end{equation}
This identity can be proved in the similar way as in Lemma 1 in
Chapter 8 of \cite{friedman-1}.  As commented in the paragraph
after Lemma 4.5 in \cite{friedman-3}, the proof of Lemma 1 can go
through when we replace Lipschitz free boundary with H\"{o}lder
continuous free boundary with the H\"{o}lder exponent
$\alpha>\frac12$.

Now we will take the derivative  of (\ref{eq:u(x,t)}) with
respect to $x$ to obtain
\begin{eqnarray}
 \frac{\partial}{\partial x} w(x,t) =\sum_{i=0}^2 \frac{\partial}{\partial x}
 M_i(x,t) \label{eq:partial_x w}
\end{eqnarray}
and let $x\downarrow b(t)$. Since $\partial_x w(b(s)+,s)$ and
$\sigma(b(s),s)$, $t_0\leq s<t$, are continuous and $b(t)$ is
H\"{o}lder continuous with exponent $\alpha>\frac12$ (see
Theorem~\ref{theorem:holder-cont}), taking $\rho(s)=\frac12
\sigma^2(b(s), s)\,\partial_x w(b(s)+,s)$ in (\ref{eq:jump
identity}), we obtain
\begin{equation}\label{eq:partial M0}
\begin{split}
\lim_{x\downarrow b(t)} \frac{\partial}{\partial x} M_0(t)
&=\lim_{x\downarrow b(t)} \frac{\partial}{\partial x} \int_{t_0}^t
ds\,  \frac12 \sigma^2(b(s),s)\,\partial_x w(b(s)+,s)G(x,t;b(s),s)\\
 &= \frac14 \sigma^2(b(t),t)\,
\partial_x w(b(t)+,t) + \int^t_{t_0} ds\, \frac12 \sigma^2(b(s),s)\,\partial_x w(b(s)+,s)
\,\partial_x G(b(t),t;b(s),s).
\end{split}
\end{equation}
On the other hand, by Lemmas \ref{lemma:lip-holder} and
\ref{lemma:ut-bounded}, $w(y,t_0)$ and $h(y,s)$ are bounded in
$\R\times [t_0,T]$. Using the Dominated Convergence Theorem we get
\begin{eqnarray}
 \lim_{x\downarrow b(t)} \frac{\partial}{\partial x}M_1(x,t) &=&
 \int_{b(t)}^{+\infty} dy\, \partial_x G(b(t),t;y,t_0) w(y,t_0)
 \triangleq N_1(t), \label{eq:def-N1}\\
 \lim_{x\downarrow b(t)} \frac{\partial}{\partial x}M_2(x,t) &=&
 \int_{t_0}^tds \int_{b(s)}^{+\infty} dy\,
 \partial_x G(b(t),t;y,s)h(y,s) \triangleq
 N_2(t), \label{eq:def-N2}
\end{eqnarray}
It follows from (\ref{eq:partial_x w}) - (\ref{eq:def-N2}) that
$\partial_x w(b(t)+,t)$ satisfies (\ref{eq:volterra}).

Let us prove the converse. For any solution $v(t)$ of the integral equation
(\ref{eq:volterra}), we can define $w(x,t)$ as follows
\begin{eqnarray}\label{eq:defnofw}
w(x,t) := - \int_{t_0}^t ds\, v(s)\, \frac12 \sigma^2(b(s),s)\,
G(x,t;b(s),s) + \int_{b(t)}^{+\infty} dy\, G(x,t;y,t_0) w(y,t_0)
+\int_{t_0}^{t}ds \int_{b(s)}^{+\infty} dy\,
 G(x,t;y,s)h(y,s), \nonumber\\
\hspace{5cm} t_0\leq t\leq T, x\geq b(t), \label{eq:def-w}
\end{eqnarray}
and $w(x,t_0) :=
\partial_t u(x,t_0)$.
We will show in the following that $w(x,t)$ is a classical
solution of (\ref{eq:w_1}) and that  $t \rightarrow \partial_x w(b(t)+,t)$ is continuous.

Now we will show that $w(x,t)$ defined in (\ref{eq:def-w}) is a
classical solution of (\ref{eq:w_1}) on $[t_0,T]$ with initial
condition $\partial_tu(x,t_0)$. By definition $w(x,t_0) =
\partial_t u(x,t_0)$. On the other hand we have that $\lim_{x\rightarrow +\infty}
w(x,t)=0$, which follows from the facts that $\lim_{x\rightarrow
+\infty}G(x,t;y,t_0) =0$ and $\sigma$, $v(s)$, $w(y,t_0)$ and
$h(y,s)$ are all bounded. Furthermore, using the properties of the
Green function and the definition of $w$ (see \ref{eq:defnofw}),
we also have that $\L_{\D} w(x,t) = h(x,t)$ for $x>b(t)$, $t\in
[t_0,T]$. Observe that  $\partial_t w$ , $\partial_x w$ and
$\partial^2_x w$ all exist and are all continuous in this domain.

In the following we will show that $\partial_x w(b(t)+,t)=v(t)$,
which implies the continuity of $\partial_x w(b(t)+,t)$. We
differentiate $w(x,t)$ with respect to $x$ and let $x\downarrow
b(t)$. Since $v(t)$ and $\sigma$ are continuous and $b(t)$ is
H\"{o}lder continuous with exponent $\alpha>\frac12$, we can apply
the jump identity (\ref{eq:jump identity}) with $\rho(s)=\frac12
\sigma^2(b(s),s) v(s)$. Following the steps that lead to
(\ref{eq:volterra}) in the first part of the proof, we obtain
\begin{equation}\label{eq:inteq-w}
 \partial_x w(b(t)+,t) =-\frac14 \sigma^2(b(t),t)\,v(t) -\int_{t_0}^t ds\, v(s) \frac12 \sigma^2(b(s),s)\, \partial_x G(b(t), t; b(s),s) + \sum_{i=1}^2 N_i(t).
\end{equation}
Comparing (\ref{eq:inteq-w}) to (\ref{eq:volterra}), we see that
$\partial_x w(b(t)+,t)= v(t)$, $t_0\leq t\leq T$.

Then it remains to show that $w(b(t),t)=0$, $t_0\leq t\leq T$. To
this end, since we have already shown $\L_{\D}w=h$, $w$ satisfies
the Green's identity given by (\ref{eq:Green-identity}).
Integrating the identity (\ref{eq:Green-identity}) and using
(\ref{eq:defnofw}) and the fact that $\lim_{x\rightarrow
+\infty}w(x,t)=0$ we can write
\begin{equation}\label{eq:w-identity-1}
\begin{split}
&\int_{t_0}^t ds\, w(b(s),s) \left[\left(\frac12 \sigma^2(b(s),s)+
\sigma \sigma_x(b(s),s)\right) \,
\partial_y G(x,t;b(s),s) - \left(\mu-\frac12 \sigma^2(b(s),s)\right)G(x,t;
b(s), s)\right]=0, \\
& \hspace{12cm} x>b(t), t_0\leq t\leq T.
\end{split}
\end{equation}
Let $x>b(t)$. Integrating both sides of (\ref{eq:w-identity-1})
on $[x,+\infty)$ and using the fact that $\partial_x G= -\partial_y G$, we
obtain
\begin{equation*}
\begin{split}
0&=\int_{t_0}^t ds\, w(b(s),s) \left[- \left(\frac12
\sigma^2(b(s),s)+ \sigma \sigma_x \right) \int_x^{+\infty}
du\,\partial_x G(u,t;b(s),s)-\left(\mu-\frac12
\sigma^2(b(s),s)\right)\int_x^{+\infty} du\,
G(u,t;b(s),s)\right]\\
&=\int_{t_0}^t ds\, w(b(s),s) \left[\left(\frac12
\sigma^2(b(s),s)+ \sigma \sigma_x\right)\, G(x,t;b(s),s) -
\left(\mu-\frac12 \sigma^2(b(s),s)\right)\,\int_x^{+\infty} du\,
G(u,t;b(s),s)\right].
\end{split}
\end{equation*}
Taking the derivative with respect to $x$, letting $x\downarrow
b(t)$ and using the jump identity (\ref{eq:jump identity}) with\\
$\rho(s)= \left(\frac12 \sigma^2(b(s),s)+ \sigma \sigma_x\right)\,
w(b(s),s)$, we arrive at
\begin{equation}\label{eq:w-identity-2}
\begin{split}
&\frac12 \left(\frac12 \sigma^2(b(s),s)+ \sigma
\sigma_x(b(s),s)\right)\, w(b(t),t) \\
&= \int_{t_0}^t ds\, w(b(s),s)\left[\left(\frac12
\sigma^2(b(s),s)+ \sigma \sigma_x(b(s),s)\right)
\partial_y G(b(t),t;b(s),s) - \left(\mu-\frac12
\sigma^2(b(s),s)\right)G(b(t),t;b(s),s) \right].
\end{split}
\end{equation}
Since $b(t)$ is H\"{o}lder continuous with exponent $\alpha>1/2$,
we have
\[
 \left|\partial_y G(b(t),t; b(s),s)\right| \leq \frac{C}{(t-s)^{\frac32
 -\alpha}}.
\]
Therefore both $\partial_y G(b(t),t;b(s),s)$ and
$G(b(t),t;b(s),s)$ are integrable. Consequently, it follows from
(\ref{eq:w-identity-1}), (\ref{eq:w-identity-2}) and the Dominated
Convergence Theorem that $w(b(t),t)=0$, $t_0\leq t \leq T$. \hfill
$\square$

\begin{proof}[Proof of Lemma~\ref{lemma:u_xt cont}]
Proof of (i). Let $v(t)$ be the unique continuous solution of the Volterra
equation (\ref{eq:volterra}). Define $w(x,t)$ as in
(\ref{eq:def-w}). The Lemma~\ref{lemma:correspondence} shows that
$w(x,t)$ is a classical solution to equation (\ref{eq:w_1}). Let
us define
\[
 \tilde{u}(x,t) = u(x,t_0) + \int_{t_0}^t w(x,s) ds, \quad x\geq
 b(t), t_0\leq t\leq T.
\]
It is easy to check that $\tilde{u}(x,t)$ is a classical solution
of the equation (\ref{eq:u_1}) - (\ref{eq:u_4}) with initial
condition $u(x,t_0)$. Since (\ref{eq:u_1}) - (\ref{eq:u_4}) has a
unique solution,  we conclude that $u(x,t) = \tilde{u}(x,t)$,
$x\geq b(t)$ and $t_0\leq t\leq T$.
 Lemma~\ref{lemma:correspondence} also implies that
\[
 \partial_x\partial_t u(b(t)+,t) = \partial_x w(b(t)+, t) = v(t),
 \quad t_0\leq t\leq T,
\]
 which implies that $\partial_x\partial_t u(b(t)+,t)$,
$t_0\leq t\leq T$, is continuous. The statement follows since
$t_0>0$ is arbitrary. \\

\noindent Proof of (ii). Let $(x,t)$ be such that $x>b(t)$. Choosing $t_0<t$
such that $b(t_0)<x$, we can see that $\int_{t_0}^t ds
\partial_x G(x,t;b(s),s)< +\infty$. As a result, we have
\[
 \frac{\partial}{\partial x}M_0(x,t) =\int^t_{t_0}
ds\,\frac12 \sigma^2(b(s),s) \,\partial_x w(b(s)+,s)
\partial_x G(x,t;b(s),s).
\]
We have shown in part (i) that $\partial_x w(b(s)+,s)$ is
continuous with respect to $s$. It is easy to show $\partial_x
M_0(x,t)$ is continuous around a sufficiently small neighborhood of
$(x,t)$. One can also show that the functions $\partial_x
M_i(x,t)$, $i \in \{1,2\}$ are also continuous by similar means. Thus, it is clear
from (\ref{eq:partial_x w}) that $\partial_x \partial_t u(x,t)$ is
continuous in this small neighborhood around $(x,t)$. Therefore, the
part (ii) of Lemma~\ref{lemma:u_xt cont} follows, because of the
arbitrary choice of $x$ and $t$.
\end{proof}

\subsection{Proof of Proposition \ref{prop:un_t} (i)}
\label{app:A.3}

We will use the following result in Lemma 4.1 in page 239 of
\cite{friedman-2}:
\begin{lemma}\label{lemma:compactness}
For any $a<b<\log{K}$, $0<t_1<t_2<T$, if both $u(x,t)$ and
$\partial_t u(x,t)$ belong to $L^2((t_1,t_2); L^2(a,b))$, then
$u(t)$ belongs to $C((t_1,t_2); L^2(a,b))$.
\end{lemma}
In this lemma, $L^2((t_1,t_2); L^2(a,b))$ is the class of $L^2$
maps which map $t\in (t_1,t_2)$ to the Hilbert space $L^2(a,b)$.
On the other hand $C((t_1,t_2); L^2(a,b))$ is the class of
continuous maps which map $t\in (t_1,t_2)$ to $L^2(a,b)$.

 The proof of (\ref{eq:u_n_t}) is similar to that of (\ref{eq:cont-at-b-us}):
 First, we will study the penalty problem associated
to the free boundary problem (\ref{eq:u_n_1}) - (\ref{eq:u_n_3}).
Then, we will list some key estimates for the solution of the
penalty problem. And finally using Lemma~\ref{lemma:compactness}
we will conclude. We will give a sketch of this proof below.

Let us consider the following penalty problem
\begin{equation}\label{eq:penalty-1}
\begin{split}
 &\L_{\D} u_n^{\epsilon} +
 \beta_{\epsilon}(u_n^{\epsilon}-g_{\epsilon})= f^{\epsilon}_n(x,t), \quad
 x\in\R, \,0<t<T, \\
 &u_n^{\epsilon}(x,0)=g_{\epsilon}(x), \quad x\in \R,
\end{split}
\end{equation}
in which $0<\epsilon<1$, $g_{\epsilon}(x)\in C^{\infty}(\R)$ such
that $g_{\epsilon}(x)= (K-e^x)^+$ when $x$ satisfies $|K-e^x|\geq
\epsilon$. We define $f_n^{\epsilon}(x,t)= \zeta_{\epsilon} \ast
f_n(x,t)$, where $\zeta_{\epsilon}$ is the standard mollifier in
$x$ and $t$ (see \cite{evans} Appendix C4 in page 629). As a
result, we have $f_n^{\epsilon}(x,t)\in C^{\infty} (\R\times
(0,T))$. Moreover, because $f_n(x,t)$ is continuous,
$f_n^{\epsilon}(x,t)$ uniformly converge to $f_n(x,t)$ on any
compact domains as $\epsilon \rightarrow 0$. On the other hand,
from our assumption that $\partial_t u_{n-1}(x,t)$ is bounded for
any $\epsilon >0$ and $\nu$ is a probability measure on $R$, we
obtain that
\begin{equation}
 \partial_t f_n(x,t) \text{ is bounded in } \R\times[\epsilon,T], \quad \text{ for any } \epsilon>0 \label{condition}.
\end{equation}
Thanks to (\ref{condition}), it is easy to see that $\partial_t
f_n^{\epsilon}(x,t)$ are uniformly bounded for any $\epsilon>0$.
The penalty functions $\beta_{\epsilon}(x)$  is a sequence of
infinitely differentiable, negative, increasing and concave
functions such that $\beta_{\epsilon}(0)=-C_{\eps} \leq
-(r+\lambda)K - r \epsilon$. The limit of the sequence is
\[
\lim_{\epsilon \rightarrow 0}\beta_{\epsilon}(x)= \begin{cases}0,
& x\geq 0,
\\ -\infty, & x<0. \end{cases}
\]

It is well known that the penalty problem has a classical solution
(see page 1009 of \cite{fri-kin}). Moreover, a proof similar to
that of the proof of Theorem 2.1 of \cite{yang} shows that
$u_n^{\epsilon}(x,t) \in C^{\infty}(\R\times(0,T))\cap
L^{\infty}(\R\times(0,T))$.

On the other hand, $u^{\epsilon}_n(x,t)$ satisfy the following
estimates for any $a<b<\log{K}$, $0<t_1<t_2 \leq T$,
\begin{eqnarray}
&&\int_a^b \left(\frac{\partial u_n^{\epsilon}}{\partial
t}\right)^2 (x,t)dx \leq C, \quad t\in
[t_1,t_2],\label{penalty-est-1}\\
&&\int_{t_1}^{t_2}\int_a^b \left(\frac{\partial^2
u_n^{\epsilon}}{\partial x\partial t}\right)^2 dx dt \leq C,
\label{penalty-est-2}\\
&&\int_{t_1}^{t_2}\int_a^b \left(\frac{\partial^2
u_n^{\epsilon}}{\partial t^2}\right)^2 dx dt + \int_a^b
\left(\frac{\partial^2 u_n^{\epsilon}}{\partial x\partial
t}\right)^2(x,t) dx \leq C, \quad t\in
[t_1,t_2],\label{penalty-est-3}
\end{eqnarray}
in which $C$ is a constant independent of $\epsilon$. These
estimates use similar techniques to the ones used in the proofs of
Lemmas 2.8, 2.10 and 2.11 in \cite{yang}, since $f_n(x,t)$
satisfies (\ref{condition}). (Similar estimates can also be found
in \cite{fri-kin}). We will give the proof for the inequality
(\ref{penalty-est-3}) below. The other inequalities can be
similarly obtained.

\begin{proof}[Proof of inequality (\ref{penalty-est-3})]
 Let us
consider $w_n(x,t)=\partial_t u_n^{\epsilon}(x,t)$. Since
$u_n^{\epsilon}(x,t)\in C^{\infty}(\R\times(0,T))$, it follows
from (\ref{eq:penalty-1}) that $w_n(x,t)$ satisfies
\begin{equation}\label{eq:penalty-wn}
 \L_{\D} w_n + \beta'_{\epsilon}(u_n^{\epsilon} -g_{\epsilon})w_n
 = \frac{\partial}{\partial t}f^{\epsilon}_n(x,t).
\end{equation}
Let $\eta(x,t)\in C_0^{\infty}(\R\times (0,T))$, such that
$\eta(x,t)=1$ for $(x,t)\in [a,b]\times[t_1,t_2]$, and
$\eta(x,t)=0$ outside a small neighborhood of
$[a,b]\times[t_1,t_2]$. Multiplying both sides of
(\ref{eq:penalty-wn}) by $\eta^2 \partial_t w_n$ and integrating
over the domain $\Omega_t=\R\times(0,t)$ in which $t_1\leq t\leq
t_2$, we obtain
\begin{eqnarray*}
0 &=&\int\int_{\Omega_t}\eta^2 \left(\frac{\partial w_n}{\partial
 t}\right)^2 dx ds - \int\int_{\Omega_t}\frac12 \sigma^2 \eta^2 \frac{\partial^2
 w_n}{\partial x^2} \frac{\partial w_n}{\partial t} dx ds -
 \int\int_{\Omega_t} \left(\mu-\frac12 \sigma^2\right) \eta^2
 \frac{\partial w_n}{\partial x} \frac{\partial w_n}{\partial t}
 dx ds \\
 && + (r+\lambda) \int\int_{\Omega_t}\eta^2
 w_n\frac{\partial w_n}{\partial t} dx ds + \int\int_{\Omega_t} \eta^2
 \beta'_{\epsilon}(u^{\epsilon}_n - g_{\epsilon})w_n \frac{\partial
 w_n}{\partial t} dx ds -\int\int_{\Omega_t}\eta^2 \frac{\partial
 w_n}{\partial t} \frac{\partial}{\partial t} f^{\epsilon}_n(x,s) dx ds \\
 &\triangleq& I_1 + I_2 + I_3 + I_4 + I_5 + I_6,
\end{eqnarray*}
where $I_j$ is the j-th term on the left and $\sigma=\sigma(x,t)$
satisfying the assumption \eqref{eq:ass-sigma}. In the following,
we will estimate each $I_j$ separately. In deriving these
estimates we will make use of the inequality
\begin{equation}\label{sim-ineq}
\frac16 A^2 +AB + \frac96 B^2 \geq 0,
\end{equation}
for any $A, B \in \mathbb{R}$. In the following estimations, $C$
will represent different constants independent of $\epsilon$.
\begin{eqnarray*}
 I_2 &=& - \frac12 \int\int_{\Omega_t} \sigma^2 \eta^2 \frac{\partial^2 w_n}{\partial
 x^2}\frac{\partial w_n}{\partial t} dx ds
 = \frac12 \int\int_{\Omega_t} \sigma^2
 \eta^2 \frac{\partial w_n}{\partial x}\frac{\partial^2 w_n}{\partial x\partial
 t} dx ds + \int\int_{\Omega_t} \sigma \eta \frac{\partial \sigma\eta}{\partial
 x}\frac{\partial w_n}{\partial x} \frac{\partial w_n}{\partial
 t}dx ds \\
 &=& \frac14 \int\int_{\Omega_t} \sigma^2 \eta^2 \frac{\partial}{\partial
 t}\left(\frac{\partial w_n}{\partial x}\right)^2 dx ds +
 \int\int_{\Omega_t}\sigma\eta \frac{\partial \sigma\eta}{\partial x}\frac{\partial
 w_n}{\partial x}\frac{\partial w_n}{\partial t} dx ds\\
 &=& \frac14 \int_{\R} \sigma^2\eta^2 \left(\frac{\partial w_n}{\partial
 x}\right)^2 (x,t) \,dx - \frac12 \int\int_{\Omega_t} \sigma\eta\frac{\partial
 \sigma\eta}{\partial t}\left(\frac{\partial w_n}{\partial x}\right)^2
 dx ds + \int\int_{\Omega_t} \sigma\eta\frac{\partial \sigma\eta}{\partial
 x}\frac{\partial w_n}{\partial x}\frac{\partial w_n}{\partial t}
 dx ds\\
 &\geq& \frac{\delta^2}{4} \int_{\R} \eta^2 \left(\frac{\partial w_n}{\partial
 x}\right)^2 (x,t) dx -  \frac12 \left|\sigma\eta\frac{\partial \sigma\eta}{\partial
 t}\right|_{L^{\infty}} \int\int_{\Omega_t} \left(\frac{\partial^2 u_n^{\epsilon}}{\partial x\partial
 t}\right)^2 dx ds \\
 && - \frac96\int\int_{\Omega_t}\left(\sigma\frac{\partial \sigma\eta}{\partial
 x}\right)^2 \left(\frac{\partial^2 u_n^{\epsilon}}{\partial x\partial
 t}\right)^2 dx ds - \frac16\int\int_{\Omega_t}
 \eta^2\left(\frac{\partial w_n}{\partial t}\right)^2 dx ds\\
 &\geq& \frac{\delta^2}{4} \int_{\R} \eta^2 \left(\frac{\partial w_n}{\partial
 x}\right)^2 (x,t) dx -C - \frac16
 \int\int_{\Omega_t}\eta^2\left(\frac{\partial w_n}{\partial
 t}\right)^2 dx ds.
\end{eqnarray*}
The first four equalities follow from integration by part. The
first inequality follows from the assumption \eqref{eq:ass-sigma}
and the inequality (\ref{sim-ineq}) with $A=\eta\frac{\partial
w_n}{\partial t}$ and $B=\sigma\frac{\partial \sigma\eta}{\partial
x}\frac{\partial w_n}{\partial x}$. The last inequality follows
from estimation (\ref{penalty-est-2}).

For $I_i$ (i=3, 4, 5), a similar procedure yields
\begin{eqnarray*}
I_3\geq -C-\frac16\int\int_{\Omega_t} \eta^2 \left(\frac{\partial
w_n}{\partial t}\right)^2 dx ds,\, I_4\geq
-C-\frac16\int\int_{\Omega_t} \eta^2 \left(\frac{\partial
w_n}{\partial t}\right)^2 dx ds,\, I_5 \geq
-C-\frac16\int\int_{\Omega_t} \eta^2 \left(\frac{\partial
w_n}{\partial t}\right)^2 dx ds.
\end{eqnarray*}
For $I_6$, we have
\begin{eqnarray*}
 I_6 &=& - \int\int_{\Omega_t} \eta^2
 \frac{\partial w_n}{\partial t} \frac{\partial}{\partial t} f_n^{\epsilon} dx ds
 \geq -\frac96\int\int_{\Omega_t} \eta^2
 \left( \frac{\partial}{\partial t} f_n^{\epsilon} \right)^2 dx ds -\frac16\int\int_{\Omega_t}
 \eta^2 \left(\frac{\partial w_n}{\partial t}\right)^2 dx ds\\
 &\geq& -C -\frac16\int\int_{\Omega_t}
 \eta^2 \left(\frac{\partial w_n}{\partial t}\right)^2 dx ds.
\end{eqnarray*}
The first inequality can be obtained using (\ref{sim-ineq}),
whereas to obtain the last inequality, we use the fact that
$\partial_t f_n^{\epsilon}(x,t)$ is uniformly bounded. Combining
all these estimates for $I_j$, we obtain
\[
\frac16\int\int_{\Omega_t} \eta^2\left(\frac{\partial
w_n}{\partial t}\right)^2 dx ds + \frac{\delta^2}{4} \int_{\R}
\eta^2 \left(\frac{\partial w_n}{\partial x}\right)^2 (x,t) dx
\leq C.
\]
This completes the proof of (\ref{penalty-est-3}). \hfill
$\square$

Using a similar proof to that of Lemma 2.2 of \cite{yang}, we can
show that $u_n^{\epsilon}(x,t)$ is uniformly bounded. Thus there
is a subsequence that $\{u_n^{\epsilon_k}\}$ converges weakly to
$u_n$ in $L^2((a,b);L^2(t_1,t_2))$ for any $a<b<\log{K}$,
$0<t_1<t<t_2<T$ (see Appendix D in \cite{evans} for an account of
the concept of weak convergence). On the other hand, it follows
from the estimates in (\ref{penalty-est-1}) -
(\ref{penalty-est-3}) that $\frac{\partial
u_n^{\epsilon}}{\partial t}$ and $\frac{\partial^2
u_n^{\epsilon}}{\partial x \partial t}$ are uniformly bounded in
$L^2(a,b)$, $\frac{\partial^2 u_n^{\epsilon}}{\partial x
\partial t}$ and $\frac{\partial^2 u_n^{\epsilon}}{\partial t^2}$
are uniformly bounded in $L^2((t_1, t_2); L^2(a,b))$. Therefore
there exists a further subsequence satisfying
\[
\frac{\partial u_n^{\epsilon_{k_j}}}{\partial t} \rightharpoonup
\frac{\partial u_n}{\partial t}, \quad
\frac{\partial^2u_n^{\epsilon_{k_j}}}{\partial x\partial
t}\rightharpoonup \frac{\partial^2 u_n}{\partial x\partial
t},\quad \frac{\partial^2 u^{\epsilon_{k_j}}_n}{\partial
t^2}\rightharpoonup \frac{\partial^2 u_n}{\partial t^2},
\]
where derivatives of $u_n$ are defined in weak sense (see Appendix
D in \cite{evans}). Here, the convergences are weak convergences.
Since  $||u|| \leq \liminf_j ||u_n^{\epsilon_{k_j}}||$ (see
Appendix D in \cite{evans} ) (\ref{penalty-est-1}) -
(\ref{penalty-est-3}) imply that
\[
 \frac{\partial u_n}{\partial t} \in L^{\infty} ((t_1,t_2);
 L^2(a,b)), \quad \frac{\partial^2 u_n}{\partial t^2} \in L^2((t_1,t_2);
 L^2(a,b)).
\]
Then it follows from Lemma \ref{lemma:compactness} that the
derivative $\partial_t u_n$ exists and is inside the space
$C((t_1,t_2); L^2(a,b))$. On the other hand, for fixed $t\in
[t_1,t_2]$, it also follows from (\ref{penalty-est-1}) and
(\ref{penalty-est-3}) and
 the Sobolev Embedding Theorem (see, for example, Theorem 4 in page 266 of \cite{evans}) that
\begin{equation}\label{eq:hold-unt}
 \left|\frac{\partial u_n}{\partial t} (x, t)-  \frac{\partial u_n}{\partial t} (\bar{x}, t)\right| \leq C |x-\bar{x}|^{1/2}, \quad x, \bar{x} \in (a,b),
\end{equation}
in which $C$ is a positive constant that does not depend on $t$.
We already know that $\partial_t u_n(\cdot, t)$ is a continuous
map with respect to $t$, therefore (\ref{eq:hold-unt}) implies
that
\[
 \frac{\partial u_n}{\partial t}  \in C((a,b)\times (t_1,t_2)).
\]
Therefore $\partial_t u_n \in C\left(\R\times (0,T]\right)$
because the choice of $a, b, t_1$ and $t_2$ are arbitrary and
$\partial_t u_n \in C\left([\log{K}, +\infty)\times (0,T]\right)$
since $[\log{K}, +\infty)\times (0,T]\in \C_n$. Moreover, we have
\begin{equation}\label{eq:lim-2}
 \lim_{x\downarrow b(t_0)}\frac{\partial u_n}{\partial t}(x,t_0) =
 \lim_{t\rightarrow t_0^-} \frac{\partial u_n}{\partial
 t}(b(t_0),t)=0,
\end{equation}
because $(b_n(t_0),t)$ is inside the stopping region for $t<t_0$
as $b_n(t)$ is decreasing.

\end{proof}

\bibliographystyle{abbrvnat}
\bibliography{regular_ref}

\begin{thebibliography}{24}
\providecommand{\natexlab}[1]{#1}
\providecommand{\url}[1]{\texttt{#1}}
\expandafter\ifx\csname urlstyle\endcsname\relax
  \providecommand{\doi}[1]{doi: #1}\else
  \providecommand{\doi}{doi: \begingroup \urlstyle{rm}\Url}\fi

\bibitem[Achdou(2008)]{achdou}
Y.~Achdou.
\newblock An inverse problem for a parabolic variational inequality with an
  integro-differential operator.
\newblock \emph{SIAM Journal on Control and Optimization}, 47\penalty0
  (2):\penalty0 733--767, 2008.

\bibitem[Bayraktar(2008)]{bayraktar-finite-horizon}
E.~Bayraktar.
\newblock A proof of the smoothness of the finite time horizon {A}merican put
  option for jump diffusions.
\newblock \emph{To appear in the SIAM Journal on Control and Optimization},
  2008.
\newblock Available at \texttt{http://arxiv.org/abs/math.OC/0703782}.

\bibitem[Cannon et~al.(1974)Cannon, Henry, and Kotlow]{Cannon}
J.~R. Cannon, D.~B. Henry, and D.~B. Kotlow.
\newblock Continuous differentiability of the free boundary for weak solutions
  of the {S}tefan problem.
\newblock \emph{Bulletin of the American Mathematical Society}, 80:\penalty0
  45--48, 1974.

\bibitem[Chen and Chadam(2006/07)]{chen-chadam}
X.~Chen and J.~Chadam.
\newblock A mathematical analysis of the optimal exercise boundary for
  {A}merican put options.
\newblock \emph{SIAM Journal on Mathematical Analysis}, 38\penalty0
  (5):\penalty0 1613--1641 (electronic), 2006/07.

\bibitem[Cont and Tankov(2004)]{cont}
R.~Cont and P.~Tankov.
\newblock \emph{Financial modelling with jump processes}.
\newblock Chapman \& Hall/CRC Financial Mathematics Series. Chapman \&
  Hall/CRC, Boca Raton, FL, 2004.

\bibitem[Evans(1998)]{evans}
L.~C. Evans.
\newblock \emph{Partial differential equations}, volume~19 of \emph{Graduate
  Studies in Mathematics}.
\newblock American Mathematical Society, Providence, RI, 1998.

\bibitem[Friedman(1964)]{friedman-1}
A.~Friedman.
\newblock \emph{Partial differential equations of parabolic type}.
\newblock Prentice-Hall Inc., Englewood Cliffs, N.J., 1964.

\bibitem[Friedman(1976)]{friedman-2}
A.~Friedman.
\newblock \emph{Stochastic differential equations and applications. {V}ol. 2}.
\newblock Academic Press [Harcourt Brace Jovanovich Publishers], New York,
  1976.
\newblock Probability and Mathematical Statistics, Vol. 28.

\bibitem[Friedman(1975)]{friedman-3}
A.~Friedman.
\newblock Parabolic variational inequalities in one space dimension and
  smoothness of the free boundary.
\newblock \emph{J. Functional Analysis}, 18:\penalty0 151--176, 1975.

\bibitem[Friedman and Kinderlehrer(1974/75)]{fri-kin}
A.~Friedman and D.~Kinderlehrer.
\newblock A one phase {S}tefan problem.
\newblock \emph{Indiana University Mathematics Journal}, 24\penalty0
  (11):\penalty0 1005--1035, 1974/75.

\bibitem[Friedman and Shen(2002)]{friedman-shen}
A.~Friedman and W.~Shen.
\newblock A variational inequality approach to financial valuation of
  retirement benefits bases on salary.
\newblock \emph{Finance and Stochastics}, 6\penalty0 (3):\penalty0 273--302,
  2002.

\bibitem[Garroni and Menaldi(1992)]{garr-mena-1}
M.~G. Garroni and J.-L. Menaldi.
\newblock \emph{Green functions for second order parabolic integro-differential
  problems}, volume 275 of \emph{Pitman Research Notes in Mathematics Series}.
\newblock Longman Scientific \& Technical, Harlow, 1992.

\bibitem[Lady{\v{z}}enskaja et~al.(1968)Lady{\v{z}}enskaja, Solonnikov, and
  Uralchva]{lad}
O.~A. Lady{\v{z}}enskaja, V.~A. Solonnikov, and N.~N. Uralchva.
\newblock \emph{Linear and Quasi-linear Equations of Parabolic Type}.
\newblock American Mathematical Society, Providence, Rhode Island, 1968.

\bibitem[Lamberton and Mikou(2008)]{lamberton}
D.~Lamberton and M.~Mikou.
\newblock The critical price for the {A}merican put in an exponential
  {L}\'{e}vy model.
\newblock \emph{Finance Stochastics}, 12:\penalty0 561--581, 2008.

\bibitem[Levendorski{\u\i}(2004)]{leven}
S.~Z. Levendorski{\u\i}.
\newblock Pricing of the {A}merican put under {L}\'evy processes.
\newblock \emph{International Journal of Theoretical and Applied Finance},
  7\penalty0 (3):\penalty0 303--335, 2004.

\bibitem[Lieberman(1996)]{lieberman}
G.~M. Lieberman.
\newblock \emph{Second order parabolic differential equations}.
\newblock World Scientific Publishing Co. Inc., River Edge, NJ, 1996.

\bibitem[Peskir(2005)]{peskir}
G.~Peskir.
\newblock On the {A}merican option problem.
\newblock \emph{Mathematical Finance}, 15\penalty0 (1):\penalty0 169--181,
  2005.

\bibitem[Peskir(2007)]{peskir-angle}
G.~Peskir.
\newblock Principle of smooth fit and diffusions with angles.
\newblock \emph{Stochastics}, 79\penalty0 (3-4):\penalty0 293--302, 2007.

\bibitem[Pham(1997)]{pham}
H.~Pham.
\newblock Optimal stopping, free boundary, and {A}merican option in a
  jump-diffusion model.
\newblock \emph{Applied Mathematics and Optimization}, 35\penalty0
  (2):\penalty0 145--164, 1997.

\bibitem[Pham(1998)]{pham-vis}
H.~Pham.
\newblock Optimal stopping of controlled jump diffusion processes: a viscosity
  solution approach.
\newblock \emph{Journal of Mathematical Systems, Estimations, and Control},
  8\penalty0 (1):\penalty0 1--27, 1998.

\bibitem[Rust(1934)]{rust}
W.~Rust.
\newblock A theorem on volterra integral equations of the second kind with
  discontinuous kernels.
\newblock \emph{American Mathematical Monthly}, 41\penalty0 (6):\penalty0
  346--350, 1934.

\bibitem[Schaeffer(1976)]{schaeffer}
D.~G. Schaeffer.
\newblock A new proof of the infinite differentiability of the free boundary in
  the {S}tefan problem.
\newblock \emph{Journal of Differential Equations}, pages 266--269, 1976.

\bibitem[Yang et~al.(2006)Yang, Jiang, and Bian]{yang}
C.~Yang, L.~Jiang, and B.~Bian.
\newblock Free boundary and {A}merican options in a jump-diffusion model.
\newblock \emph{European Journal of Applied Mathematics}, 17\penalty0
  (1):\penalty0 95--127, 2006.

\bibitem[Zhang(1997)]{zhang}
X.~L. Zhang.
\newblock Numerical analysis of {A}merican option pricing in a jump-diffusion
  model.
\newblock \emph{Mathematics of Operations Research}, 22:\penalty0 668--690,
  1997.

\end{thebibliography}
\end{document}